\documentclass[11pt]{article}
\usepackage{amssymb}
\usepackage{amsfonts}
\usepackage{color}
\usepackage{xy}
\xyoption{all}
\usepackage{setspace}

\newtheorem{theorem}{Theorem}[section]
\newtheorem{proposition}[theorem]{Proposition}
\newtheorem{lemma}[theorem]{Lemma}
\newtheorem{corollary}[theorem]{Corollary}
\newtheorem{definition}[theorem]{Definition}

\newtheorem{example}[theorem]{Example}

\newtheorem{remark}[theorem]{Remark}

\newenvironment{proof}{{\noindent\bf Proof.}}{\hfill $\Box$\par\vskip3mm}

\newcommand{\Ker}{{\rm Ker}\,}

\newcommand{\im}{{\rm Im}\,}

\newcommand{\Hom}{{\rm Hom}}
\newcommand{\End}{{\rm End}}

\newcommand{\Cc}{\mathcal{C}}
\newcommand{\Dd}{\mathcal{D}}
\newcommand{\Ee}{\mathcal{E}}

\newcommand{\Mm}{\mathcal{M}}

\newcommand{\Tt}{\mathcal{T}}
\newcommand{\Ss}{\mathcal{S}}
\def\RR{{\mathbb R}}

\def\NN{{\mathbb N}}
\def\CC{{\mathbb C}}
\def\ZZ{{\mathbb Z}}

\begin{document}
\baselineskip16pt
\author{Miodrag Cristian Iovanov \small \\ University of Southern California, 
Department of Mathematics, \\ 3620 South Vermont Ave. KAP 108 \\ Los Angeles, California 90089-2532\\ and \\  Department of Mathematics, University of Bucharest \\ Str. Academiei 14, Bucharest 1, 010014 Romania }
\title{Abstract integrals in algebra} 

\maketitle
	
\begin{abstract}
We generalize the results on existence and uniqueness of integrals from compact groups and Hopf algebras in a pure (co)algebraic setting, and find a series of new results on (quasi)-co-Frobenius and semiperfect coalgebras. For a coalgebra $C$, we introduce the generalized space of integrals $\int_M=\Hom^C(C,M)$ associated to a right $C$-comodule $M$ and study connections between ``uniqueness of integrals'' $\dim(\int_M)\leq \dim(M)$ and ``existence of integrals'' $\dim(\int_M)\geq \dim(M)$ for all $M$ and representation theoretic properties of $C$: (quasi)-co-Frobenius, semiperfect. We show that a coalgebra is co-Frobenius if and only if existence and uniqueness of integrals holds for any finite dimensional $M$. We give the interpretation for $\int_M$ for the coalgebra of representative functions of a compact group - they will be "quantum"-invariant vector integrals. As applications, new proofs of well known characterizations of co-Frobenius coalgebras and Hopf algebras are obtained, as well as the uniqueness of integrals in Hopf algebras. We also give the consequences for the representation theory of infinite dimensional algebras. We give an extensive class of examples which show that the results of the paper are the best possible. These examples are then used to give all the previously unknown connections between the various important classes of coalgebras appearing in literature.
\end{abstract}

\maketitle

\section*{Introduction}

Let $G$ be a compact group. It is well known that there is a unique (up to multiplication) left invariant Haar measure $\mu$ on $G$, and a unique left invariant Haar integral. If $H$ is a Hopf algebra over a field $K$, an element $\lambda\in H^*$ is called a left integral for $H$ if $\alpha\lambda=\alpha(1)\lambda$ for all $\alpha\in H^*$. For a compact group $G$, let $R_c(G)$ be the $\CC$-coalgebra (Hopf algebra) of representative functions on $G$, consisting of all $f:G\rightarrow \CC$ such that there are (continuous) $u_i,v_i:G\rightarrow \CC,\,i=\overline{1,n}$ such that $f(xy)=\sum\limits_{i=1}^nu_i(x)v_i(y)$ for all $x,y\in G$. Then $\int$ restricted to $R_c(G)$ is an integral in the Hopf algebra (coalgebra) sense (see for example \cite[Chapter 5]{DNR}). The uniqueness of integrals for compact groups has a generalization for Hopf algebras: if a nonzero (left) integral exists in $H$, then it is shown by Radford \cite{R} that it is unique, in the sense that the dimension of the space of left integrals equals $1$. 

\vspace{.5cm}

For a Hopf algebra $H$, it is easy to see that a left integral $\lambda$ is the same as a morphism of right $H$-comodules (left $H^*$-modules) from $H$ to the right $H$-comodule $K$ with comultiplication $K\ni a\mapsto a\otimes 1_H\in K\otimes H$. Then it is natural to generalize this definition to arbitrary finite dimensional $H$-comodules by putting $\int_M=\Hom_{H^*}(H,M)$. The advantage of this definition is that it can be considered for arbitrary coalgebras, since in general, for a coalgebra $C$ there is no canonical $C$-comodule structure on $K$. We give an explicit description of the space of these generalized integrals for the case of the representative coalgebra (Hopf algebra) of a (locally) compact group and an interpretation at the group level. More precisely, we will consider vector-valued integrals $\int$ on $G$, $\int: C(G)\rightarrow \CC^n=V$ (or $\int:L^1(G)\rightarrow \CC^n$) with the "quantum-invariance" property $\int x\cdot f=\eta(x)\cdot \int f$ for any $f\in R_c(G)$, where $\eta:G\rightarrow \End(V)=\End(\CC^n)$; it turns out that $\eta$ must actually be a representation of $G$. Then $V$ with the left $G$-action is turned naturally into a right $R_c(G)$-comodule and the integral restricted to $R_c(G)$ turns out to be an algebraic integral in the above sense, that is, $\int\in\Hom^{R_c(G)}(R_c(G),V)$.

\vspace{.5cm}

We note that in the case of a locally compact group $G$, the coalgebra structure of $R_c(G)$ is the one encoding the information of representative functions, and so of the group itself: the comultiplication of $R_c(G)$ is defined by $\Delta(f)=\sum\limits_{i=1}^nu_i\otimes v_i$, for the above $u_i,v_i$ such that $f(xy)=\sum\limits_{i=1}^nu_i(x)v_i(y)$ for all $x,y\in G$. The algebra structure is given by $(f*g)(x)=f(x)g(x)$, and this comes by "dualizing" the comultiplication $\delta$ of the coalgebra structure of $\CC[G]$, which is defined as $\delta(x)=x\otimes x$ for $x\in G$. Since this coalgebra structure does not involve the group structure of $G$ in any way ($G$ might as well be a set), it is to be expected that only the coalgebra structure of $R_c[G]$ will encapsulate information on $G$. This suggests that a generalization of the existence and uniqueness of integrals results should be possible for the case of coalgebras.

\vspace{.5cm}

With this in mind, we generalize the existence and uniqueness results from Hopf algebras to the pure coalgebraic setting. For a coalgebra $C$ and a right finite dimensional right $C$-comodule $M$ we define the space of left integrals $\int_{l,M}=\Hom_{C^*}(C^C,M^C)$ and similarly for left $C$-comodules $N$ let $\int_{r,N}=\Hom_{C^*}({}^CC,{}^CN)$ be the space of right integrals. We note that this definition has been considered before in literature; see \cite[Chapter 5.4]{DNR}. It is noted there that if $C$ is a left and right co-Frobenius coalgebra, then $\dim(\Hom_{C^*}(C,M))\leq \dim(\Hom_{C^*}(C,M))$; this result was proved in \cite{St} for certain classes of co-Frobenius coalgebras (finite dimensional, or cosemisimple, or which are Hopf algebras). Our goal is to prove here far more general results, and give the generalization of the now well known result of Hopf algebras stating that a Hopf algebra is co-Frobenius if and only if it has nonzero left integrals (equivalently, has right integrals), and in this case, the integral is unique up to scalar multiplication.

\vspace{.5cm}

It is natural to think to the dimensional comparison $\dim(\int_{l,M})\leq\dim(M)$ as a "uniqueness" of integrals for $M$ and then to the statement $\dim(\int_{l,M})\geq \dim(M)$ as "existence of integrals". Recall that a coalgebra is called left (right) co-Frobenius if $C$ embeds in $C^*$ as left (right) $C^*$-modules, and simply co-Frobenius if it is both left and right co-Frobenius. We first show that for a coalgebra which is (just) left co-Frobenius, the "uniqueness" of (left) integrals holds for all right $C$-comodules $M$ ($\dim(\int_{l,M})\leq\dim(M)$) and the "existence" of (right) integrals holds as well for all left $C$-comodules $N$ ($\dim(\int_{r,N})\leq\dim(N)$). Examples are provided later on to show that the converse statements do not hold (even if both left and right existence - or both left and right uniqueness - of integrals are assumed). On the way, we produce some interesting characterizations of the more general quasi-co-Frobenius (shortly, QcF) coalgebras; these will show that the co-Frobenius and quasi-co-Frobenius properties are fundamentally about a certain duality between the left and right indecomposable components of $C$, and the multiplicities of these in $C$ (Propositions \ref{p.1} and \ref{p.2}). 

\vspace{.5cm}

One main result of the paper is the Theorem \ref{TheTh}, which extends the results from Hopf algebras; it states that a coalgebra is left and right co-Frobenius if and only if existence and uniqueness of left integrals hold for all (right) $C$-comodules $M$ and equivalently, for all left comodules $M$. This adds to the previously known symmetric characterization of co-Frobenius coalgebras from \cite{I}, where it is shown that $C$ is co-Frobenius if and only if $C$ is isomorphic to its left (or, equivalently, to its right) rational dual $Rat({}_{C^*}C^*)$. Moreover, it is shown there that this is further equivalent to the functors $C^*$-dual $\Hom(-,C^*)$ and $K$-dual $\Hom(-,K)$ from ${}_{C^*}\Mm$ to ${}_K\Mm$ being isomorphic when restricted to the category of left (equivalently, right) rational $C^*$-modules which is the same as that of right $C$-comodules: $Rat({}_{C^*}\Mm)=\Mm^C$. This brings up an interesting comparison to the algebra case: if the two functors were to be isomorphic on the whole category of left $C^*$-modules, one would have that $C^*$ is a Frobenius algebra (by well known facts of Frobenius algebras, see \cite{CR}), so $C^*$ (and $C$) would be finite dimensional. This showed why the co-Frobenius coalgebra concept is a generalization of the Frobenius algebra in the infinite dimensional case. Here, the above mentioned Theorem \ref{TheTh} allows us to extent this view by giving a new interesting characterization of co-Frobenius coalgebras: $C$ is co-Frobenius if and only if the functors $C^*$-dual $\Hom(-,C^*)$ and $K$-dual $\Hom(-,K)$ are isomorphic (only) on the subcategory of ${}_{C^*}\Mm$ consisting of finite dimensional rational left $C^*$-modules (and then further equivalent to the right hand side version of this statement). In fact, quite interestingly, we note that for $C$ to be co-Frobenius, it is enough for these two functors to be isomorphic when evaluated in vector spaces, but by an isomorphism which is not necessarily natural; the existence of a natural isomorphism follows thereafter (Theorem \ref{thCatcF}). We also give a categorical characterization of co-Frobenius coalgebras which generalizes the fact that a finite dimensional $K$-algebra ${}_A\Mm$ is Frobenius if and only of the forgetful functor ${}_A\Mm\rightarrow {}_K\Mm$ is Frobenius, i.e. has isomorphic left and right adjoints (Corollary \ref{ThAdjFr}). Another aplication is the fact that $C$ is QcF if and only if the classes of projective and injective right (or equivalently, left) comodules coincide (\ref{ThPrIn}), which is analogous to another characterization of (finite dimensional) QF algebras. As further applications, we find the well known equivalent characterizations of Hopf algebras with nonzero integrals of Lin, Larson, Sweedler, Sullivan as well as the uniqueness of integrals as a consequence of our general results.

We also give an extensive class of examples which will show that all the results in the paper are the best possible (Section 3). On the side, we also obtain interesting examples (of one sided and two sided) semiperfect, QcF and co-Frobenius coalgebras showing that all possible inclusions between these classes are strict (for example, we note that there are left and right QcF coalgebras which are left co-Frobenius but not right co-Frobenius, or which are neither left nor right co-Frobenius). In particular, it is known that a left QcF is left semiperfect, and we prove a new and interesting fact: a left QcF coalgebra is also right semiperfect (Theorem \ref{qcfsp}). Our examples of coalgebras are associated to graphs and are usually subcoalgebras of the path coalgebra. We also find a similar functorial characterization of semiperfect coalgebras: $C$ is left semiperfect if and only if the forgetful functor from finite dimensional left $C$-comodules to $K$-vector spaces is the restriction of a representable functor ${}^C\Mm\rightarrow {}_K\Mm$. Integrals for algebras and consequences for the representation theory of infinite dimensional algebras (Section 4).

In Section 5, we look at the abstract spaces of integrals in the case of the representative Hopf algebra (coalgebra) of compact groups $G$, and note that the abstract integrals are in fact restrictions of unique vector integrals $\int$ on $C(G)$ - the algebra of complex continuous functions on $G$ - which have a certain "quantum"-invariance: $\int(x\cdot f)dh=\eta(x)\int(f)dh$, where $\eta$ is a finite dimensional representation of $G$. In particular, we note a nice short Hopf algebra proof of a well known fact (due to Peter and Weyl) stating that any finite dimensional representation of a compact group is completely reducible, and give the statements on the existence and uniqueness of "quantum" integrals for compact groups.

\section{The General results}

For basic facts on coalgebras and their comodules we refer the reader to \cite{A, DNR, M, S}. Recall that if $C$ is a coalgebra, then $C=\bigoplus\limits_{S\in \Ss}E(S)^{n(S)}$ as left $C$-comodules, where $\Ss$ is a set of representatives of left $C$-comodules, $E(S)$ is an injective envelope of the left comodule $S$ contained in $C$ and $n(T)$ are natural (positive) numbers. Similarly, $C=\bigoplus\limits_{T\in\Tt}E(T)^{p(T)}$ as right comodules, with $0\neq p(T)\in\NN$ and $\Tt$ a set of representatives for the right simple $C$-comodules. We convey to use the letter $S$ whenever a left simple comodule is inferred and $T$ for the right simple $C$-comodules. We use $\Mm^C$ (and ${}^CM$ respectively) for the category of right (respectively left) $C$-comodules, and ${}_{A}\Mm$ or $\Mm_A$ for the categories of left or right $A$-modules over a ring $A$.

Also, we write $M^C$ or ${}_{C^*}M$ whenever we refer to the structure of $M$ as a right $C$-comodule or of the structure of $M$ as a left $C^*$-module. If $M\in\Mm^C$ then $M$ has a left $C^*$-module structure defined by $c^*\cdot m=c^*(m_1)m_0$, where  for the comultiplication $\rho:M\rightarrow M\otimes C$ of $M$ we use the Sweedler notation with the summation symbol omitted: $\rho(m)=m_0\otimes m_1\in M\otimes C$. We always have in mind this identification of the right $C$-comodules with the so called rational left $C^*$-modules which is in fact an equivalence of categories $\Mm^C\simeq Rat-{}_{C^*}\Mm$. If $M\in M^C$ (so $M\in {}_{C^*}\Mm$) then $M^*$ has a natural structure of a right $C^*$-module $(m^*\cdot c^*)(m)=m^*(c^*\cdot m)=m^*(m_0)c^*(m_1)$, which will be always understood when talking about $M^*$ as a right $C^*$-module. Moreover, if $M$ is a finite dimensional comodule, $M^*$ is a rational (finite dimensional) right $C^*$-module so it has a compatible left $C$-comodule structure (i.e. $M^*\in {}^C\Mm$). The left $C^*$-module (or right $C$-comodule) $N^*$ for $N\in {}^C\Mm$ (or finite dimensional $N\in {}^C\Mm$) can be defined similarly.

\begin{definition}\label{def}
Let $M$ be a right $C$-comodule. The space of the left integrals of $M$ will be $\int_{l,M}=\Hom(C^C,M^C)$, the set of morphisms or right $C$-comodules (left $C^*$-modules), regarded as a left $C^*$-module by the action $(c^*\cdot \lambda)(c)=\lambda(c\cdot c^*)=\lambda(c^*(c_1)c_2)$. Similarly, if $N\in {}^CM$ then $\int_{r,N}=\Hom({}^CC,{}^CN)$ is a right $C^*$-module. 
\end{definition}

We will sometimes write just $\int_M$ or $\int_N$ if there is no danger of confusion, that is, if the comodule $M$ or $N$ has only one comodule structure (for example, it is not a bimodule). The following Lemma is very useful for understanding these spaces of integrals and for computations, and the proof is standard:

\begin{lemma}\label{l.1}
If $M\in \Mm^C$ and $N\in {}^C\Mm$ then $\Hom_{C^*}(M,N^*)\simeq \Hom_{C^*}(N,M^*)$ naturally in $M$ and $N$; more precisely $\Hom({}_{C^*}M,{}_{C^*}N^*)\simeq \Hom(N_{C^*},M^*_{C^*})$.
\end{lemma}
\begin{proof}
This follow from the usual Hom-Tensor adjunctions: $$\Hom_{C^*}(M,N^*)\simeq \Hom_K(N\otimes_{C^*}M,K) \simeq\Hom_{C^*}(N,M^*)$$
\end{proof}

By the above, since for a finite dimensional comodule $M$ we have $M\simeq M^*{}^*$, we have that $\int_{l,M}=\Hom(C^C,M^C)\simeq\Hom(C,M^*{}^*)=\Hom(M^*,C^*)=\Hom(M^*,Rat(C^*_{C^*}))$.\\
We recall the following definitions here, and refer the reader to \cite{L} or \cite{DNR} for more details and characterizations of these coalgebras. 

\begin{definition}
Let $C$ be a coalgebra.
\begin{itemize}
\item
$C$ is called left (right) co-Frobenius if $C$ embeds in $C^*$ as left (right) $C^*$-modules. 
\item
$C$ is called left (right) QcF (quasi-co-Frobenius) if $C$ embeds in a coproduct of copies of $C^*$ as left (right) $C^*$-modules, $C\hookrightarrow \bigoplus\limits_IC^*$. 
\item
$C$ is called co-Frobenius (or QcF) if it is both left and right co-Frobenius (or left and right QcF). 
\item
A coalgebra $C$ is said to be left semiperfect if $E(T)$ is finite dimensional for all $T\in\Tt$ (injective envelopes of simple right comodules are finite dimensional); right semiperfect coalgebras are defined similarly. 
\end{itemize}
\end{definition}

We have the following two characterizations of QcF and co-Frobenius coalgebras which show that these two properties are actually representation theoretic properties of a coalgebra, and they refer to a certain duality between the indecomposable left and the indecomposable right components of $C$ and their multiplicities in $C$.

\begin{proposition}\label{p.1}
Let $C$ be a coalgebra. Then $C$ is left QcF if and only if for all $T\in\Tt$ there is some (unique) $S\in \Ss$ such that $E(T)\simeq E(S)^*$.
\end{proposition}
\begin{proof}
First note that given $T\in\Tt$, such $S\in \Ss$ is unique: if $E(T)\simeq E(S)^*$ then $E(T)$ is rational and finitely generated so $E(T)$ is finite dimensional and therefore $E(S)\simeq (E(S)^*)^*\simeq E(T)^*$ so $S$ is the socle of $E(T)^*$. Moreover, if this happens, denoting by $\Ss'$ is the set of those $S$'s with $E(S)^*\simeq E(T)$ for some $T$, we have morphisms of $C^*$-modules 
\begin{eqnarray*}
C=\bigoplus\limits_{T\in\Tt}E(T)^{p(T)} & \hookrightarrow & \bigoplus\limits_{T\in\Tt}\bigoplus\limits_{\NN}E(T)\simeq \bigoplus\limits_\NN\bigoplus\limits_{S\in\Ss'}E(S)^*\simeq \bigoplus\limits_\NN\bigoplus\limits_{S\in\Ss'}E(S)^*{}^{n(S)}\\
& \hookrightarrow & \bigoplus\limits_\NN\prod\limits_{S\in\Ss}E(S)^*{}^{n(S)} =\coprod\limits_\NN C^*
\end{eqnarray*}
since $C^*\simeq \prod\limits_{S\in\Ss}E(S)^*{}^{n(S)}$ as left $C^*$-modules. Conversely, if $C$ is left QcF, by definition, we have a monomorphism $C\stackrel{\varphi}{\hookrightarrow}\prod\limits_{i\in L}E(S_i)^*$, with $S_i$ simple left comodules and $E(S_i)$ their injective envelopes (because $C^*\simeq \prod\limits_{S\in\Ss}E(S)^*{}^{n(S)}$). 
Since $E(T)$ are finite dimensional (see \cite[Chapter 3.2-3.3]{DNR}) and they are direct summands of $C$, it is straightforward to see that this monomorphism restricts to a monomorphism $E(T)\hookrightarrow \prod\limits_{i\in I}E(S_i)^*$ with $I$ a finite subset of $L$. To see this, let $p_i$ be the projection onto $E(S_i)^*$; then $\Ker(\varphi\vert_{E(T)})=\bigcap\limits_{i\in L}\Ker(p_i\circ\varphi\vert_{E(T)})=0$, so there is a finite $I\subset L$ such that $\bigcap\limits_{i\in I}\Ker(p_i\circ\varphi\vert_{E(T)})=0$ since $E(T)$ is finite dimensional, yielding such a monomorphism. Now $E(T)$ are also injective as $C^*$-modules (see \cite[Chapter 2.4]{DNR} for example) and thus this monomorphism splits: $E(T)\oplus X\simeq \bigoplus\limits_{i\in I}E(S_i)^*$. This shows that $E(T)$ is projective (which is also known by the fact that $C$ is left QcF). By \cite[Lemma 1.4]{I}, $E(S_i)^*$ are local indecomposable cyclic with unique maximal submodule $M_i$, with $E(S_i)^*/M_i\simeq S_i^*$. If $J(M)$ denotes the Jacobson radical of the module $M$, then we have 
$$\frac{E(T)\oplus X}{J(E(T)\oplus X)}=\frac{E(T)}{J(E(T))}\oplus \frac{X}{J(X)}=\frac{\oplus_{i\in I}E(S_i)^*}{\oplus_{i\in I}M_i}=\bigoplus\limits_{i\in I}\frac{E(S_i)^*}{M_i}=\bigoplus\limits_{i\in I}S_i^*$$
But $E(T)/J(E(T))\neq 0$ since $E(T)$ is finite dimensional, and then the composition $f=\,(E(T)\rightarrow \frac{E(T)}{J(E(T))}\hookrightarrow \bigoplus\limits_{k}S_k^*\rightarrow S_i^*)$ is nonzero for at least one $i$ (all the morphisms are the natural ones). Then the diagram:
$$
\xymatrix{
& E(T)\ar[dl]_{\overline{f}}\ar[d]^f & \\
E(S_i)^* \ar[r]_p & S_i^* \ar[r] & 0 
}
$$
is completed commutatively by $\overline{f}$ which has to be surjective (otherwise ${\rm Im}(\overline{f})\subseteq M_i$ so $f=p\overline{f}=0$). But $E(S_i)^*$ is projective, so $\overline{f}$ must split and we get $E(T)\simeq E(S_i)^*\oplus Y$ ($Y\subseteq E(T)$), and then $Y=0$ since $E(T)$ is indecomposable. This completes the proof.
\end{proof}

\begin{proposition}\label{p.2}
Let $C$ be a left QcF coalgebra and $\sigma:\Tt\rightarrow \Ss$ be defined such that $\sigma(T)=S$ if and only if $E(T)\simeq E(S)^*$ (this is well defined by the above Proposition). Then $C$ is left co-Frobenius if and only if $n(\sigma(T))\geq p(T)$, $\forall T\in\Tt$.
\end{proposition}
\begin{proof}
If $C\hookrightarrow C^*$ in ${}_{C^*}\Mm$ then for each $T\in\Tt$ there is a monomorphism $\varphi:E(T)^{p(T)}\hookrightarrow\prod\limits_{S\in\Ss}E(S)^*{}^{n(S)}$. 
$C$ is semiperfect since it is left QcF (see \cite{NT1}) and so the $E(T)$'s are finite dimensional. Therefore, as before, we may again find a finite subset $F$ of $\Ss$ and a monomorphism $E(T)^{p(T)}\hookrightarrow \bigoplus\limits_{S\in F}E(S)^*$ which splits as in the proof of Proposition \ref{p.1}: 
$$E(T)^{p(T)}\oplus Y=\bigoplus\limits_{S\in F}E(S)^*{}^{n(S)}$$ 
Again, since all the $E(S)^*$ are local cyclic indecomposable, we get that $E(T)\simeq E(S)^*$ for some $S\in F$; moreover, there have to be at least $p(T)$ indecomposable components isomorphic to $E(T)$ on the right hand side of the above equation. But since $E(S)^*$ and $E(S')^*$ are not isomorphic when $S$ and $S'$ are not, we conclude that we must have $n(S)\geq p(T)$ for the $S$ for which $E(S)^*\simeq E(T)$, i.e. $S=\sigma(T)$ and $n(\sigma(T))\geq p(T)$.\\
Conversely, if $p(T)\leq n(\sigma(T))$ we have monomorphisms of left $C^*$-modules
$$C=\bigoplus\limits_{T\in\Tt}E(T)^{p(T)}\hookrightarrow \bigoplus\limits_{T\in\Tt}E(T)^{n(\sigma(T))}\hookrightarrow \bigoplus\limits_{S\in\Ss}E(S)^*{}^{n(S)}\subseteq \prod\limits_{S\in\Ss}E(S)^*{}^{n(S)}=C^*$$
\end{proof}

Let $C_S=\sum\limits_{S'\simeq S, S'\leq C}S'$ be the simple subcoalgebra of $C$ associated to $S$. Then $C_S$ is a simple coalgebra which is finite dimensional, and $C_S\simeq S^{n(S)}$. The dual algebra $C_S^*$ of $C_S$ is a simple finite dimensional algebra, $C_S^*=(S^*)^{n(S)}$ as left $C_S^*$-modules (or $C^*$-modules) and thus $C_S^*\simeq M_{n(S)}(\Delta_S)$, where $\Delta_S=\End_{C^*}(S^*)$ is a division algebra. By Lemma \ref{l.1} we also have $\Delta_S\simeq \End(S_{C^*})$, and it is easy to see that the isomorphism in Lemma \ref{l.1} also preserves the multiplicative structure thus giving an isomorphism of algebras. Let $d(S)=\dim(\Delta_S)$. Then, as $C_S^*\simeq M_{n(S)}(\Delta_S)=(S^*)^{n(S)}$, we have $\dim(C_S)=\dim(C_S^*)=d(S)\cdot n(S)^2=n(S)\cdot \dim S$ and therefore $\dim(S)=\dim(S^*)=n(S)d(S)$. For a right simple comodule $T$ denote $d'(T)=\dim(\End({}_{C^*}T))$; note that $d'(T)=d(T^*)$ since $\End({}_{C^*}T)\simeq \End(T^*_{C^*})$ by the same Lemma \ref{l.1}. Similarly for right simple comodules $T$, $\dim(T)=d'(T)p(T)$. Then we also have $p(T)=n(T^*)$. Denote by $C_0$ the coradical of $C$; then we have that $C_0=\bigoplus\limits_{S\in\Ss}C_S$. 

\begin{remark}\label{1.rem}
Let $C$ be a left QcF coalgebra, and assume that $\End(S)=K$ for all simple left (equivalently, right) comodules $S$ (for example, this is true if $C$ is pointed or the basefield $K$ is algebraically closed). Then $C$ is left co-Frobenius if and only if $\dim(soc(E))\leq \dim(cosoc(E))$ for any finite dimensional indecomposable injective right comodule $E$, where $cosoc(E)$ represents the cosocle of $E$. Indeed, in this case, $d(S)=1=d'(T)$ and if $E(T)\simeq E(S)^*$, then $S^*=cosoc(E(T))$, so $n(\sigma(T))=n(S)=\dim(S)=\dim(cosoc(E(T)))$ and $p(T)=\dim(T)=\dim(soc(E(T)))$.
\end{remark}

\begin{proposition}\label{p.3}
The following assertions are equivalent:\\
(i) $\dim(\int_{l,M})\leq \dim(M)$ for all finite dimensional $M\in\Mm^C$.\\
(ii) $\dim(\int_{l,T})\leq \dim(T)$ for all simple comodules $T\in\Mm^C$.\\
If $C$ is a left QcF coalgebra, then these are further equivalent to \\
(iii) C is left co-Frobenius.\\
Moreover, if $C$ is left QcF then:\\
(a) $\int_{l,T}\neq 0$ for all $T\in\Tt$ if and only if $C$ is also right QcF.\\
(b) $\dim(\int_{l,T})\geq\dim(T)$ if and only if $C$ is also right co-Frobenius.
\end{proposition}
\begin{proof}
(ii)$\Rightarrow$(i) We prove (i) by induction on the length of $M$ (or on $\dim(M)$). For simple modules it holds by assumption (ii). Assume the statement holds for comodules of length less than ${\rm length}(M)$. Let $M'$ be a proper subcomodule of $M$ and $M''=M/M'$; we have an exact sequence $0\rightarrow \Hom(C,M')\rightarrow \Hom(C,M)\rightarrow \Hom(C,M'')$ and therefore $\dim(\int_{l,M})=\dim(\Hom(C^C,M))\leq\dim(\int_{l,M'})+\dim(\int_{l,M''})\leq\dim(M')+\dim(M'')$ by the induction hypothesis and thus $\dim(\int_{l,M})\leq\dim(M')+\dim(M'')=\dim(M)$.\\
(i)$\Rightarrow$(ii) is obvious.\\
Assume $C$ is left QcF and again, let $\sigma:\Tt\rightarrow \Ss$ be such that $E(T)\simeq E(\sigma(T))^*$ as given by Proposition \ref{p.1}.\\
(i)$\Leftrightarrow$(iii) Let $T_0\in\Mm^C$ be simple. Then there exists at most one $T\in\Tt$ such that $\Hom(E(T),T_0)\neq 0$. Indeed, for any $T\in\Tt$, $E(T)\simeq E(S)^*$ for $S=\sigma(T)$. Since $T_0^*$ is a rational $C^*$-module, applying Lemma \ref{l.1} we get $\Hom(E(T),T_0)=\Hom(E(T),T_0^*{}^*)=\Hom(T_0^*,E(T)^*)=\Hom(T_0^*,E(S))=\Hom(T_0^*,S)$ which is nonzero if and only if $T_0^*\simeq S=\sigma(T)$. This can only happen for at most one $T$. Thus we get that $\int_{l,T_0}=\Hom(C^C,T_0)=\Hom(\bigoplus\limits_{T\in\Tt}E(T)^{p(T)},T_0)=\prod\limits_{T\in\Tt}\Hom(E(T),T_0)^{p(T)}$ is 0 if $T_0^*$ does not belong to the image of $\sigma$, or $\int_{l,T_0}=\Hom(E(T),T_0)^{p(T)}=\Hom(T_0^*,S)^{p(T)}$ with $S=\sigma(T)=T_0^*$ as above. In this latter case, we have 
$$\dim(\int_{l,T_0})=p(T)\dim(\Hom(T_0^*,T_0^*))=p(T)d(T_0^*)$$
while $\dim(T_0)=\dim(T_0^*)=n(T_0^*)d(T_0^*)=n(\sigma(T))\cdot d(T_0^*)$. Since for $T_0^*\notin{\rm Im}(\sigma)$ $\dim(\int_{l,T_0})=0\leq \dim(T_0)$, we get that $\dim(\int_{l,T_0})\leq\dim(T_0)$ holds for all $T_0$ if and only if this takes place for $T_0$ ranging in the image of $\sigma$, and by the above equalities, this is further equivalent to $p(T)\leq n(\sigma(T))$, $\forall T\in\Tt$. By Proposition \ref{p.2} this is equivalent to $C$ being left co-Frobenius. This finishes (i)$\Leftrightarrow$(iii) under the supplementary hypothesis of $C$ being left QcF.\\
For (a) if $C$ is left QcF, since $\int_{l,T_0}\neq 0$ if and only if $T_0^*\in{\rm Im}(\sigma)$ we see that $\sigma$ is bijective if and only if $\int_{l,T}\neq 0,\,\forall T\in\Tt$ (since $\sigma$ is automatically injective). The surjectivity of $\sigma$ means that for all $S\in\Ss$, there is some $T$ such that $E(S)^*\simeq E(T)$, or $E(S)\simeq E(T)^*$, which is equivalent to $C$ being right QcF by Proposition \ref{p.1}.\\
(b) Using (a) and the above facts, $\sigma$ is bijective and so $\dim(\int_{l,T_0})=p(T)d(T_0^*)\geq \dim(T_0)=n(\sigma(T))d(T_0^*)$, ($T_0^*=\sigma(T)$) for all $T_0$ is equivalent to $n(\sigma(T))\leq p(T)$ for all $T$, that is, $p(\sigma^{-1}(S))\geq n(S)$ for all $S\in\Ss$ which means $C$ is right co-Frobenius by the right hand side version of Proposition \ref{p.2}.
\end{proof}

\begin{corollary}\label{c}
If $C$ is a left co-Frobenius coalgebra, then $\dim(\int_{l,M})\leq \dim(M)$ for all finite dimensional $M\in\Mm^C$.
\end{corollary}

\begin{remark}
The above corollary was also proved in \cite{DNT2}, under the supplimentary assumption that $C$ is also right semiperfect. In view of our later Theorem \ref{qcfsp}, this condition can, in fact, be dropped. \\
We see that by the above characterization of left QcF coalgebras, if $C$ is left QcF, then $\int_{r,S}\neq 0$ for all $S\in \Ss$; indeed, let $T=S^*$ and $S_0\in\Ss$ such that $E(T)\simeq E(S_0)^*$. Then the monomorphism $T\hookrightarrow E(S_0)^*$ produces a nonzero epimorphism $E(S_0)^*\rightarrow T^*=S\rightarrow 0$ so $\Hom(C,S)\neq 0$. Therefore, $\int_{r,N}\neq 0$ for all $N$, because any comodule $N$ contains some simple comodule $S\in\Ss$. We thus observe the following interesting
\end{remark}

\begin{corollary}
The following are equivalent for a coalgebra $C$:\\
(i) $C$ is left QcF and $\int_{l,T}\neq 0$ for all simple left rational $C^*$-modules $T$.\\
(ii) $C$ is right QcF and $\int_{r,S}\neq 0$ for all simple right rational $C^*$-modules $S$.
\end{corollary}

\begin{proposition}\label{p.4}
Let $C$ be a left co-Frobenius coalgebra. Then $\dim(\int_{r,N})\geq \dim(N)$ for all finite dimensional $N\in{}^C\Mm$. 
\end{proposition}
\begin{proof}
$C$ is also left QcF, so there is $\sigma:\Tt\rightarrow \Ss$ such that $E(T)\simeq E(\sigma(T))^*$ as in Proposition \ref{p.1}. Let $\Ss'=\sigma(\Tt)$ and $H=\bigoplus\limits_{S\in\Ss'}E(S)^{n(S)}=\bigoplus\limits_{T\in\Tt}E(T)^*{}^{n(\sigma(T))}$. Note that $H$ is projective in $\Mm_{C^*}$, since $E(T)^*$ is a projective right $C^*$-module for all $T$ (direct summand in $C^*$). Also, $C=H\oplus H'$, $H'=\bigoplus\limits_{S\in\Ss\setminus \Ss'}E(S)^{n(S)}$ and so $\dim(\Hom({}^CC,N))=\dim(\Hom_{C^*}(H,N))+\dim(\Hom_{C^*}(H',N))\geq \dim(\Hom_{C^*}(H,N))$. When $N=S_0$ a simple left $C$-comodule, then there exists exactly one $S\in\Ss'$ such that $\Hom_{C^*}(E(S),S_0)\neq 0$. Indeed, for $S\in\Ss'$, $E(S)\simeq E(T)^*$ with $S=\sigma(T)$, and since $S_0^*$ is a simple rational left $C^*$-module, 
using also Lemma \ref{l.1} we have $\Hom(E(S),S_0)=\Hom(E(S),S_0^*{}^*)=\Hom(S_0^*,E(S)^*)=\Hom(S_0^*, E(T))=\Hom(S_0^*,T)$ and this is nonzero if and only if $S_0^*\simeq T$, i.e. $S=\sigma(T_0^*)$ (all the $\Hom$ represent morphisms of $C^*$-modules). This shows that 
$$\Hom(H,S_0)=\prod\limits_{S\in\Ss'}\Hom(E(S),S_0)^{n(S)}=\prod\limits_{T\in\Tt}\Hom(S_0^*,T)^{n(\sigma(T))}=\Hom(S_0^*,S_0^*)^{n(\sigma(S_0^*))}$$
and therefore $\dim(\Hom(H,S_0))=\dim(\Hom(S_0^*,S_0^*))\cdot n(\sigma(S_0^*))=d'(S_0^*)n(\sigma(S_0^*))\geq d'(S_0^*)p(S_0^*)$ because $C$ is left co-Frobenius. But $d'(S_0^*)p(S_0^*)=\dim(S_0^*)=\dim(S_0)$ and thus we get $\dim(\Hom(H,S_0))\geq \dim(S_0)$. Since $H$ is a projective right $C^*$-module, this inequality can be extended to all finite dimensional left $C$-comodules by an inductive argument on the length of the left $C$-comodule $N$, just as in the proof of Proposition \ref{p.3}. Finally, $\dim(\int_{r,N})=\dim(\Hom(C,N))\geq \dim(\Hom(H,N))\geq \dim(N)$ and the proof is finished.\\
\end{proof}

\section{Co-Frobenius coalgebras and Hopf algebras}

The next theorem generalizes the existence and uniqueness of left and of right integrals from co-Frobenius Hopf algebras to the general case of co-Frobenius coalgebras, showing that, as in the Hopf algebra case, these are actually equivalent to the coalgebra being co-Frobenius. It is noted in \cite[Remark 5.4.3]{DNR} that for co-Frobenius coalgebras $\dim_{C^*}(C,M)\leq \dim(M)$. This was shown above in Proposition \ref{p.3} to hold in a more general case of only left co-Frobenius coalgebras (with actual equivalent conditions) and the following gives the mentioned generalization:

\begin{theorem}\label{TheTh}
A coalgebra $C$ is co-Frobenius (both on the left and on the right) if and only if $\dim(\int_{l,M})=\dim(M)$ for all finite dimensional right $C$-comodules $M$, equivalently, $\dim(\int_{r,N})=\dim(N)$ for all finite dimensional $N$ in ${}^C\Mm$.
\end{theorem}
\begin{proof}
"$\Rightarrow$" Since $C$ is left co-Frobenius, Proposition \ref{p.3} shows that $\dim(\int_{l,M})\leq \dim(M)$ for finite dimensional right comodules $M$ and as $C$ is also right co-Frobenius, the right hand side of Proposition \ref{p.4} shows that $\dim(\int_{r,M})\geq \dim(M)$ for such $M$.\\
"$\Leftarrow$" Let $T$ be a simple right $C$-comodule and $S=T^*$. Let $X$ be the socle of $Rat(C^*_{C^*})$ and $X_S=\sum\limits_{S'<C^*,S'\simeq S}S'$ be the sum of all simple sub(co)modules of $C^*$ isomorphic to $S$. It is easy to see that $X=\bigoplus\limits_{S\in\Ss}X_S$ and $X_S$ is semisimple isomorphic to a direct sum of comodules isomorphic to $S$, $X_S\simeq S^{(I)}=\coprod\limits_IS$. Then $\Hom(C,T)=\Hom(C,T^*{}^*)=\Hom(T^*,C^*)=\Hom(S,C^*)=\Hom(S,X_S)$ so $\dim(\Hom_{C^*}(C,T))=\dim(\Hom(S,X_S))$; if $I$ has cardinality greater than $n(S)$ then $\dim(\Hom(S,X_S))> \dim(\Hom(S,S^{n(S)}))=d(S)n(S)=\dim(S)=\dim(T)$ so $\dim(\Hom(C,T))>\dim(T)$ and this contradicts the hypothesis. Then we get that $I$ is finite and $\dim(\Hom(C,T))=|I|\cdot\dim(\Hom(S,S))=d(S)\cdot|I|=\dim(T)=\dim(S)=d(S)\cdot n(S)$ and thus $|I|=n(S)$. This shows that $X_S\simeq S^{n(S)}\simeq C_S$. Hence $X=\bigoplus\limits_{S\in\Ss}X_S\simeq \bigoplus\limits_{S\in\Ss}C_S\simeq C_0$ as left $C$-comodules (right $C^*$-modules).\\
Next, we show that $Rat(C^*_{C^*})$ is injective: let $0\rightarrow N'\stackrel{f}{\rightarrow} N\stackrel{g}{\rightarrow} N''\rightarrow 0$ be an exact sequence of finite dimensional left $C$-comodules; it yields the exact sequence of vector spaces $0\rightarrow \Hom(N'',Rat(C^*_{C^*}))\stackrel{g^*}{\rightarrow} \Hom(N,Rat(C^*_{C^*}))\stackrel{f^*}{\rightarrow} \Hom(N',Rat(C^*_{C^*}))$. Evaluating dimensions we get 
\begin{eqnarray*}
\dim(\Hom(N',Rat(C^*_{C^*}))) & = & \dim(\Hom(N',C^*))=\dim(\int_{l,(N')^*})=\dim(N')^*=\dim(N') \\ 
& = & \dim(N)-\dim(N'')=\dim\int_{l,N^*}-\dim\int_{l,(N'')^*}\\ 
& = & \dim\Hom(N,C^*)-\dim\Hom(N'',C^*)=\\ 
& = & \dim\Hom(N,Rat(C^*_{C^*}))-\dim\Hom(N'',Rat(C^*_{C^*}))\\
& = & \dim({\rm Im}f^*)
\end{eqnarray*}
and this shows that $f^*$ is surjective. Then, by \cite[Theorem 2.4.17]{DNR} we get that $Rat(C^*_{C^*})$ is injective. \\
Then, since $X$ is the socle of the injective left comodule $Rat(C^*_{C^*})$, we get $Rat(C^*_{C^*})\simeq E(X)$ because $X$ is essential in $Rat(C^*_{C^*})$; but $X\simeq C_0$ in ${}^C\Mm$ so $E(X)\simeq E(C_0)\simeq C$, i.e. $C\simeq Rat(C^*_{C^*})$. By \cite[Theorem 2.8]{I} we get that $C$ is left and right co-Frobenius.
\end{proof}

We now give the applications of these general results to the equivalent characterizations of co-Frobenius Hopf algebras and the existence and uniqueness of integrals for Hopf algebras. Recall that if $H$ is a Hopf algebra over a field $K$, $\lambda\in H^*$ is called a left integral for $H$ if $h^*\cdot\lambda=h^*(1)\lambda$; this is equivalent to saying that the 1-dimensional vector space $K\lambda$ is a left ideal of $H^*$ which is rational, and its right comultiplication $\rho:K\lambda\rightarrow K\lambda\otimes H$ writes $\rho(\lambda)=\lambda\otimes 1$. Let $\int_l$ denote the space of all left integrals of $H$, and defined similarly, let $\int_r$ be the space of all right integrals. Note that $\int_l=\Hom({}_{H^*}K\cdot 1,{}_{H^*}H^*)=\int_{l,K\cdot 1}$ where $K\cdot 1$ is the right $H$-comodule with comultiplication given by $1\mapsto 1\otimes 1_H$; indeed $\varphi:K\cdot 1\rightarrow H^*$, $\varphi(1)=\lambda\in H^*$, is a morphism of left $H^*$-modules if and only if $\lambda$ is an integral: $\varphi(h^*\cdot 1)=h^*\cdot \varphi(1)\Leftrightarrow h^*(1)\varphi(1)=h^*\cdot\varphi(1)$.

We will need to use the isomorphism of right $H$-comodules $\int_l\otimes H\simeq Rat({}_{H^*}H^*)$ from \cite{Sw1}, pp.330-331 (see also \cite[Chapter 5]{DNR}), which is in fact an isomorphism of $H$-Hopf modules, but we only need the comodule isomorphism (we will not use the right $H$-module structure of $Rat({}_{H^*}H^*)$). The above mentioned isomorphism is a direct easy consequence of the fundamental theorem of Hopf modules.

We note that only part of the results of the previous section (Proposition \ref{p.3} or Corollary \ref{c}) are already enough to derive the well known uniqueness of integrals for Hopf algebras.

\begin{corollary}[Uniqueness of Integrals of Hopf algebras]
Let $H$ be a Hopf algebra. Then $\dim(\int_l)\leq 1$.
\end{corollary}
\begin{proof}
If $\int_l\neq 0$, then there is a monomorphism of left $H^*$-modules $H\hookrightarrow \int_l\otimes H\simeq Rat({}_{H^*}H^*)\hookrightarrow H^*$. Therefore $H$ is left co-Frobenius and Corollary \ref{c} (or Proposition \ref{p.3}) shows that $\dim(\int_l)=\dim(\int_{l,K})\leq \dim(K)=1$.
\end{proof}

We can however derive the following more general results due to Lin, Larson, Sweedler, Sullivan \cite{L,LS,Su}.

\begin{theorem}
Let $H$ be a Hopf algebra. Then the following assertions are equivalent:\\
(i) $H$ is a left co-Frobenius coalgebra.\\
(ii) $H$ is a left QcF coalgebra.\\
(iii) $H$ is a left semiperfect coalgebra.\\
(iv) $Rat({}_{H^*}H^*)\neq 0$.\\
(v) $\int_l\neq 0$.\\
(v') $\int_{l,M}\neq 0$ for some finite dimensional right $H$-comodule $M$.\\
(vi) $\dim\int_l=1$.\\
(vii) The right hand side versions of (i)-(vi)
\end{theorem}
\begin{proof}
(i)$\Rightarrow$(ii)$\Rightarrow$(iii)$\Rightarrow$(iv) are properties of coalgebras (\cite{L}, \cite{NT1}, \cite[Chapter 3]{DNR}), (vi)$\Rightarrow$(v) is trivial and (iv)$\Leftrightarrow$(v) follows by the isomorphism $\int_l\otimes H\simeq Rat({}_{H^*}H^*)$; also (v)$\Rightarrow$(v') and (v') implies $\int_{l,M}\simeq\Hom_{H^*}(M^*,H^*)\neq 0$ so $Rat({}_{H^*}H^*)\neq 0$ ($M^*$ is rational), and thus (iv) and (v) follow. Now assume (v) holds; then (i) follows since the isomorphism of right $H$-comodules $\int_l\otimes H\simeq Rat({}_{H^*}H^*)$ shows that $H\hookrightarrow H^*$ in ${}_{H^*}\Mm^*$, i.e. $H$ is left co-Frobenius. Moreover, in this case, since $\int_r=\int_{r,K1}$, Proposition \ref{p.4} shows that $\dim(\int_r)\geq \dim(K1)=1$. In turn, by the right hand side of equivalences of (i)-(v), $H$ is also right co-Frobenius and Proposition \ref{p.3} shows that $\dim(\int_l)\leq 1$ so $\dim(\int_l)=1$ and similarly $\dim(\int_r)=1$. Hence, (v)$\Rightarrow$(i), (vi) \& (vii), and this ends the proof.
\end{proof}

\subsection*{Further applications to categorical characterizations of co-Frobenius coalgebras}

We use the above results to give a new characterization of co-Frobenius coalgebras. It is known from \cite{I} and \cite{II} that a $K$-coalgebra $C$ is co-Frobenius if and only if the functors $\Hom_{C^*}(-,C^*)$ and $\Hom_K(-,K)$ from ${}_{C^*}\Mm$ to ${}_K\Mm_{}$ are isomorphic when restricted to the category of right $C$-comodules $\Mm^C$ (rational left $C^*$-modules); see \cite{II} for similar characterizations of the more general quasi-co-Frobenius coalgebras. We will show that it is enough for the two functors to be isomorphic only on the finite dimensional rational $C^*$-modules, or even more generally, that the $C^*$-dual and the $K$-dual of any finite dimensional rational comodule have the same dimension. 

\begin{theorem}\label{thCatcF}
The following are equivalent for a coalgebra $C$:\\
(i) $C$ is co-Frobenius.\\
(ii) The functors $\Hom_{C^*}(-,C^*)$ and $\Hom_K(-,K)$ are {\it naturally} isomorphic when restricted to the category of finite dimensional right $C$-comodules.\\
(iii) The integral functor $\int_{r,-}$ on the category of finite dimensional left $C$-comodules $f.d.{}^C\Mm$ is naturally isomorphic to the forgetful functor $U:f.d.{}^C\Mm\rightarrow {}_K\Mm$.\\
(iv) The $C^*$-dual and the $K$-dual of any rational finite dimensional left $C^*$-module $M$ are isomorphic as vector spaces, that is, they have the same dimension $\dim(\Hom_{C^*}(M,C^*))=\dim(M)$.\\
(v) The right hand side version of (i)-(iii).
\end{theorem}
\begin{proof}
(i)$\Rightarrow$(ii) follows by \cite[Theorem 2.8]{I}. \\ 
(ii)$\Leftrightarrow$(iii) We have natural isomorphisms $\int_{r,N}=\Hom_{C^*}(C,N)=\Hom_{C^*}(C,(N^*)^*)\cong\Hom_{C^*}(N^*,C^*)$ and $N\cong(N^*)^*=\Hom_K(N^*,K)$ so $\int_{r,N}\cong N$ naturally for all finite dimensional left $C$-comodules $N$ if and only if $\Hom_{C^*}(N^*,C^*)\cong\Hom_K(N^*,K)$, equivalently, $\Hom_{C^*}(M,C^*)\cong\Hom_K(M,K)$ naturally for finite dimensional $M\in\Mm^C$.\\
(iii)$\Rightarrow$(iv) is obvious.\\
(iv)$\Rightarrow$(i) If $N$ is a left $C$-comodule, then $M=N^*$ is a right $C$-comodule and $\dim(M^*)=\dim(\Hom_{C^*}(M,C^*))$ and therefore 
\begin{eqnarray*}
\dim(\int_{r,N}) & = & \dim(\Hom_{C^*}({}^CC,{}^CN))=\dim(\Hom_{C^*}(C,M^*))\\
& = & \dim(\Hom({}_{C^*}M,{}_{C^*}C^*)) \,\,\,\,\,\,\, ({\rm by\,Lemma\,}\ref{l.1})\\
& = & \dim(M) \,\,\,\,\,\,\,\, ({\rm by\,hypothesis})\\
& = & \dim(N)
\end{eqnarray*}
and the result follows now as an application of Theorem \ref{TheTh}.
\end{proof}

We note how the results on integrals allow us to prove a ``Frobenius'' functor characterization for co-Frobenius coalgebras, which has eluded the theory so far. Recall that a $K$-algebra $A$ is Frobenius if and only if the forgetful functor ${}_A\Mm\rightarrow {}_K\Mm$ is Frobenius (or a strong adjunction), meaning it has isomorphic left and right adjoints. Such a statement for the forgetful functor $\Mm^C\rightarrow {}_K\Mm$ for a coalgebra $C$ automatically leads to $C$ being finite dimensional co-Frobenius (so the dual of a Frobenius algebra), so it does not recover the full generality of the co-Frobenius coalgebras. However, relaxing the strong adjunction property to a ``local'' type of property in a way similar to \cite{IV}, we obtain an interesting characterization of co-Frobenius coalgebras. We note that the results (and definitions) of \cite{IV} give characterizations of one-sided (quasi)-co-Frobenius coalgebras (corings). Recall that the forgetful functor $F:\Mm^C\rightarrow {}_K\Mm$ has a right adjoint $-\otimes C:{}_K\Mm\rightarrow \Mm^C$. Let us introduce de the

\begin{definition}
Let $(F:\Cc\rightarrow \Dd,G:\Dd\rightarrow \Cc)$ be a pair of adjoint functors between additive categories such that $\Cc$ has coproducts and $\Dd$ has products. Let us call $(F,G)$ a co-Frobenius pair, or a semi-strong adjunction if $(G,F)$ is an adjoint pair when the functors are ``restricted'' to the subcategories ${\rm fcg}(\Dd)$ and ${\rm fg(\Cc)}$ of finitely cogenerated and, respectively, finitely generated objects of $\Dd$ and, respectively, $\Cc$. More precisely, there is a natural isomorphism
$$\Hom_{\Cc}(G(Y),X)\cong \Hom_\Dd(Y,G(X))$$
for all $X$ finitely generated in $\Cc$ and $Y$ finitely cogenerated in $\Dd$.\\
(recall that $X$ finitely generated means that $X=\sum\limits_{i\in I}X_i$ implies $X=\sum\limits_{i\in F}X_i$ for some finite $F\subseteq I$, and dually, $Y$ finitely cogenerated means $\bigcap\limits_{i\in I}Y_i=0$ for subobjects of $Y$ implies $\bigcap\limits_{i\in F}Y_i=0$ for some finite $F\subseteq I$).
\end{definition}

With this we have the following nice extension of the characterization of (finite dimensional) Frobenius (co)algebras:

\begin{corollary}\label{ThAdjFr}
The following are equivalent for a coalgebra $C$:\\
(i)$C$ is co-Frobenius;\\
(ii) The functors $F:{}^C\Mm\rightarrow {}_K\Mm$ and its right adjoint $G=-\otimes C:{}_K\Mm\rightarrow {}^C\Mm$ are a co-Frobenius pair (or a semi-strong adjunction).\\
(iii) The left comodule version of (ii).
\end{corollary}
\begin{proof}
The condition in (ii) reads $\Hom^C(V\otimes C,N)\cong \Hom_K(V,N)$ for $N$ finite dimensional left comodule and $V$ finite dimensional vector space. Hence, if (ii) holds, applying it for $V=k$ it follows that $\int_{r,N}=\Hom_{C^*}(C,N)\cong \Hom_K(K,N)=N$ for all finite dimensional $N\in{}^C\Mm$ so $C$ is co-Frobenius by Theorem \ref{TheTh}. Conversely, for finite dimensional $V\in{}_K\Mm$ and $N\in{}^C\Mm$ we have a sequence of natural isomorphisms
\begin{eqnarray*}
\Hom^C(V\otimes C,N) & \cong & V\otimes \Hom^C(C,N) \cong V\otimes \Hom_{C^*}(C,(N^*)^*) \\
& \cong & V \otimes \Hom_{C^*}(N^*,C^*) \cong V\otimes \Hom_K(N^*,K)\,\,\,\,\,({\rm by\,Theorem\,}\ref{thCatcF})\\
& \cong & V \otimes \Hom_K(K,N) \cong \Hom_K(V,N)
\end{eqnarray*}
\end{proof}

\section{Examples and Applications}

We provide some examples to show that most of the general results given above are in some sense the best possible for coalgebras.

\begin{example}\label{ex1}
Let $C$ be the $K$-coalgebra having $g,c_n$, $n\geq 1, n\in\NN$ as a basis and comultiplication $\Delta$ and counit $\varepsilon$ given by
\begin{eqnarray*}
\Delta(g)=g\otimes g & , & \Delta(c_n)=g\otimes c_n+c_n\otimes g \,\,\,\forall n\\
\varepsilon(g)=1 & , & \varepsilon(c_n)=0 \,\,\,\forall n
\end{eqnarray*}
i.e. $g$ is a grouplike element and the $c_n$'s are $(g,g)$-primitive elements. Then $S=Kg$ is essential in $C$, $S$ is the only type of simple $C$-comodule and $C/S\simeq \bigoplus\limits_n(Kc_n+S)/S\simeq \bigoplus\limits_\NN S$. Then $\Hom(C/S,S)\simeq\Hom(\bigoplus\limits_\NN S,S)\simeq \prod\limits_\NN \Hom(S,S)$ and since there is a monomorphism $\Hom(C/S,S)\rightarrow \Hom(C,S)=\int_S$, it follows that $\int_S$ is infinite dimensional. (In fact, it can be seen that $\int_S\simeq \Hom(C/S,S)$: we have an exact sequence $0\rightarrow \Hom(C/S,S)\rightarrow \Hom(C,S)\rightarrow \Hom(S,S)$. The last morphism in this sequence is $0$, because otherwise it would be surjective ($\dim(\Hom(S,S))=1$) and this would imply that the inclusion $S\subseteq C$ splits, which is not the case.) 
Thus we have $\dim(S)\leq \dim\int_S$ for (all) the simple comodule(s) $S$. Then, for any $C$-comodule $N$, there exists a monomorphism $S\hookrightarrow N$ which produces a monomorphism $\int_S\hookrightarrow \int_N$ and therefore $\int_N$ is always infinite dimensional and so $\dim(\int_N)\geq \dim(N)$, for all finite dimensional $C$-comodules $N$ (and since $C$ is cocommutative, this holds on both sides). Nevertheless, $C$ is not co-Frobenius, since it is not even semiperfect: $C=E(S)$. This shows that the converse of Proposition \ref{p.4} does not hold. 
\end{example}

\begin{example}
Let $C$ be the divided power series $K$-coalgebra with basis $c_n$, $n\geq 0$ and comultiplication $\Delta(c_n)=\sum\limits_{i+j=n}c_i\otimes c_j$ and counit $\varepsilon(c_n)=\delta_{0n}$. Then $C^*\simeq K[[X]]$ - the ring of formal power series, and the only proper subcomodules of $C$ are $C_n=\bigoplus\limits_{i=0}^nKc_n$. Since all these are finite dimensional, $C$ has no proper subcomodules of finite codimension, and we have $\int_N=\Hom(C,N)=0$ for any finite dimensional $C$-comodule $N$ (again this holds both on left and on the right). Thus "uniqueness" $\dim(\int_{N})\leq \dim(N)$ holds for all $N$'s, but $C$ is not co-Frobenius since it is not even semiperfect ($C=E(Kc_0)$).
\end{example}

We give a construction which will be used in a series of examples, and will be used to show that the Propositions in the first sections are the best possible results. Let $\Gamma$ be a directed graph, with the set of vertices $V$ and the set of edges $\Ee$. For each vertex $v\in V$, let us denote $L(v)$ the set of edges coming into $v$ and by $R(v)$ the set of edges going out of $v$. For each side $m$ we denote $l(v)$ its starting vertex and $r(m)$ its end vertex: ${}^{l(v)}\bullet \stackrel{m}{\longrightarrow}\bullet^{r(v)}$. We define the coalgebra structure $K[\Gamma]$ over a field by $K$ defining $K[\Gamma]$ to be the vector space with basis $V\sqcup \Ee$ and comultiplication $\Delta$ and counit $\varepsilon$ defined by 
$$\Delta(v)=v\otimes v,\, \varepsilon(v)=1 \,{\rm for}\, v\in V;$$
$$\Delta(m)=l(m)\otimes m+m\otimes r(m)$$
Denote by $<x,y,\dots>$ the $K$-vector space with spanned by $\{x,y,\dots\}$ We note that this is the second term in the coradical filtration in the path coalgebra associated to $\Gamma$, and it is not difficult to see that this actually defines a coalgebra structure. Notice that the socle of $K[\Gamma]$ is $\bigoplus\limits_{v\in V}<v>$ and the types of simple (left, and also right) comodules are $\{<v>\mid v\in V\}$. We also note that there exist direct sum decomposition of $K[\Gamma]$ into indecomposable injective left $K[\Gamma]$-comodules 
$${}^{K[\Gamma]}K[\Gamma]=\bigoplus\limits_{v\in V}{}^{K[\Gamma]}<v;m\mid m\in L(v)>$$
and a direct sum decomposition into indecomposable right $K[\Gamma]$-comodules
$$K[\Gamma]^{K[\Gamma]}=\bigoplus\limits_{v\in V}<v;m\mid m\in R(v)>^{K[\Gamma]}$$
To see this, note that each of the components in the above decompositions is a left (respectively right) subcomodule of $K[G]$ and that it has essential socle given by the simple (left and right) $K[\Gamma]$-comodule $<v>$. For $v\in V$ let $E_r(v)=<v;m\mid m\in R(v)>^{K[\Gamma]}$ and $E_l(v)={}^{K[\Gamma]}<v;m\mid m\in L(v)>$. We have an exact sequence of right $K[\Gamma]$-comodules
$$0\rightarrow <v>^{K[\Gamma]}\rightarrow E_r(v)^{K[\Gamma]}\rightarrow \bigoplus\limits_{m\in R(v)}<r(m)>^{K[\Gamma]}\rightarrow 0$$
Since $<v>$ is the socle of $E_r(v)$, this shows that a simple right comodule $<w>$ ($w\in V$) is a quotient of an injective indecomposable component $E_r(v)$ whenever $w=r(m)$ for some $m\in R(v)$. This can happen exactly when $E_r(v)$ contains some $m\in L(w)$. Therefore we have $\Hom^{K[\Gamma]}(E_r(v),<w>)=0$ whenever $m\notin L(w)$ for any $m\in R(v)$, and $\Hom^{K[\Gamma]}(E_r(v),<w>)=\prod\limits_{m\in L(w)\mid l(m)=v}K$. Thus 
\begin{eqnarray*}
\Hom^{K[\Gamma]}(K[\Gamma],<w>) & = & \Hom^{K[\Gamma]}(\bigoplus\limits_{v\in V}E_r(v),<w>)=\prod\limits_{v\in V}\Hom^{K[\Gamma]}(E_r(v),<w>) \\ & = & \prod\limits_{v\in V}\prod\limits_{m\in L(w)\mid l(m)=v}K = \prod\limits_{m\in L(w)}<w>
\end{eqnarray*}
i.e. $\dim(\int_{l,<w>})=\dim K^{L(w)}$. Similarly, we can see that $\dim(\int_{r,<w>})=\dim K^{R(w)}$. \\
We will use this to study different existence and uniqueness of integrals properties for such coalgebras. Also, we note a fact that will be easy to use in regards to "the existence of integrals" for a coalgebra $C$: if $\dim(\int_{r,S})=\infty$ for all simple left $C$-comodules $S$, then for any finite dimensional left $C$-comodule $N$, let $S$ be a simple subcomodule of $N$; then the exact sequence $0\rightarrow \Hom^C(C,S)\rightarrow \Hom^C(C,N)$ shows that $\dim(\int_{r,N})=\infty$, so existence of right integrals holds trivially in this case.\\
We also note that the above coalgebra has the following:\\
$\bullet$ uniqueness of left (right) integrals if $|L(w)|\leq 1$ ($|R(w)|\leq 1$) for all $w\in V$, since in this case, $\dim(\int_{l,<w>})\leq 1$ for all simple right comodules $T=<w>$, and this follows by Proposition \ref{p.3}\\
$\bullet$ existence of left (right) integrals if $|L(w)|=\infty$ ($|R(w)|=\infty$) for all $w\in V$ since then $\dim(\int_{l,w})=\dim(K^{L(w)})=\infty$ and it follows from above.\\
$\bullet$ $K[\Gamma]$ is left (right) semiperfect if and only if $R(w)$ ($L(w)$) is finite for all $w\in V$ (if $R(w)$ ($L(w)$) is infinite for some $w\in V$ then $K[\Gamma]$ is not left (right) semiperfect since $E_r(v)$ is not finite dimensional- it contains the elements of $R(v)$ in a basis). Therefore, when this fails, $K[\Gamma]$ cannot be left (right) QcF nor left (right) co-Frobenius. \\
$\bullet$ If $|R(w)|\geq 2$ for some $w\in V$, then $K[\Gamma]$ is not left QcF. Otherwise, $E_r(w)\simeq E_l(v)^*$, with $E_l(v)=<v;m\mid m\in L(v)>$ with both $E_r(w), E_l(v)$ finite dimensional; but $E_l(v)$ has socle $<v>$ of dimension 1, so $E_r(w)\simeq E_l(v)^*$ is local by duality. But $\dim(E_r(w)/<w>)=|R(w)|\geq 2$ and $E_r(w)/<w>$ is semisimple, so it has more than one maximal subcomodule, which is a contradiction. Similarly, if $|L(w)|\geq 2$ then $K[\Gamma]$ is not right QcF (nor co-Frobenius).

\begin{example}\label{ex2.5}
Let $\Gamma$ be the graph
$$
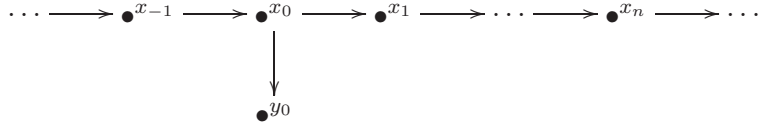
\xymatrix{
\dots \ar[r] & \bullet^{x_{-1}} \ar[r] & \bullet^{x_0} \ar[r] \ar[d] & \bullet^{x_1} \ar[r] & \dots \ar[r] & \bullet^{x_n} \ar[r] & \dots \\
& & \bullet^{y_0} & & & &
}$$
and $C=K[\Gamma]$. By the above considerations, we see that $C$ has the existence and uniqueness property of left integrals of simple modules: $\dim(\int_{l,T})=\dim(T)=1$ for all right simple $C$-comodules $T$. But this coalgebra is not left QcF (nor co-Frobenius) because $|R(x_0)|=2$ and it is also not right QcF, because $E_l(y_0)$ is not isomorphic to a dual of a right injective $E_r(v)$, as it can be seen directly by formulas, or by noting that $E_r(x_0)^*=<x_0,[x_0x_1],[x_0y_0]>^*$ and $E_r(y_0)=<y_0>^*$ are the only duals of right injective indecomposables containing the simple left comodule $<y_0>$, and they have dimensions 3 and 1 respectively.\\
This shows that the characterization of co-Frobenius coalgebras from Theorem \ref{TheTh} cannot be extended further to requiring existence and uniqueness only for simple comodules, as it in the case of Hopf algebras, where existence for the simple comodule $K1$ is enough to infer the co-Frobenius property.
\end{example}

\begin{example}\label{ex3}
Consider the poset $V=\bigsqcup\limits_{n\geq 0}\NN^n$ with the order given by the "levels" diagram
$${0}\rightarrow \NN\rightarrow \NN\times\NN\rightarrow \NN\otimes\NN\otimes\NN\rightarrow \dots$$
and for elements in consecutive levels we have that two elements are comparable only in the situation $(x_0,x_1,\dots,x_n)<(x_0,x_1,\dots,x_n,x)$, $x_0=0$, $x_1,\dots,x_n,x\in\NN$. This is making $V$ into a poset which is actually a tree with root $v_0=(0)$. Visually, we can see this as in the diagram (the arrows indicate ascension):
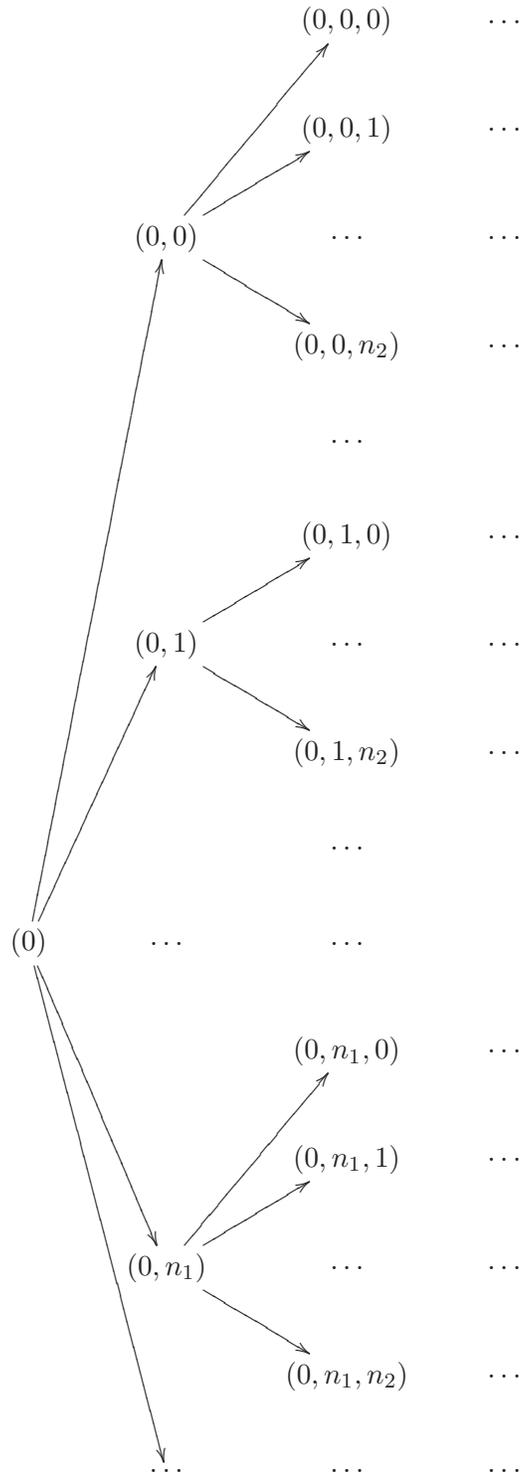
{
$$
\xymatrix{
& & (0,0,0) & \dots \\
& & (0,0,1) & \dots \\ 
& (0,0) \ar[uur]\ar[ur]\ar[dr] & \dots & \dots\\
& & (0,0,n_2) & \dots \\
& & \dots & \\
& & (0,1,0) & \dots \\
& (0,1) \ar[ur]\ar[dr] & \dots & \dots\\
& & (0,1,n_2) & \dots \\
& & \dots & \\
(0)\ar[uuuuuuur]\ar[uuur]\ar[dddr]\ar[dddddr] & \dots & \dots \\ 
& & (0,n_1,0) & \dots \\
& & (0,n_1,1) & \dots \\ 
& (0,n_1) \ar[uur]\ar[ur]\ar[dr] & \dots & \dots\\
& & (0,n_1,n_2) & \dots \\
& \dots & \dots & \dots
}
$$
}
Let $\Gamma$ be the above tree, i.e. having vertices $V$ and sides (with orientation) given by two consecutive elements of $V$. For each pair of consecutive vertices $a,b$ we have exactly one side $[ab]$ and the comultiplication reads \\
$\bullet$ $\Delta(a)=a\otimes a$ and $\varepsilon(a)=1$ for $a\in V$ (i.e. $a$ is a grouplike element)\\
$\bullet$ $\Delta([ab])=a\otimes [ab]+[ab]\otimes b$ and $\varepsilon([ab])=0$ for $[ab]\in M$, that is, $b\in S(a)$ (i.e. $[ab]$ is ($a,b$)-primitive)\\
We see that here we have $|L(v)|\leq 1$ for all $v\in V$ (in fact $|L(v)|=1$ for $v\neq v_0$ and $|L(v_0)|=0$) so left uniqueness of integrals holds: $\dim\int_{l,M}\leq \dim M$, for all finite dimensional rational left $K[\Gamma]^*$-modules $M$ by Example \ref{ex2.5}. Since $|R(v)|=\infty,\,\forall v\in V$, by the same \ref{ex2.5} it follows that $\dim(\int_{r,N})\geq\dim N$ for $N\in{}^{K[\Gamma]}\Mm$ (existence of right integrals). However, this coalgebra is not left co-Frobenius (nor QcF) because $|R(v)|=\infty,\,\forall v\in V$. This shows that the converse of Proposition \ref{p.4} and Corollary \ref{c} combined does not hold. More generally, for this purpose, we could consider an infinite rooted tree (that is, a tree with a pre-chosen root) with the property that each vertex has infinite degree.\\
We also note that this coalgebra is not right co-Frobenius (nor QcF) either, because the dual of a left injective indecomposable comodule cannot be isomorphic to a right injective indecomposable comodule, since the latter are all infinite dimensional.
\end{example}

\begin{example}
As seen in the previous example, it is also not the case that "left uniqueness" and "right existence" of integrals imply the fact that $C$ is right co-Frobenius; this can also be seen because there are coalgebras $C$ that are left co-Frobenius and not right co-Frobenius (see \cite{L} or \cite[Chapter 3.3]{DNR}). Then the left existence and right uniqueness hold by the results in Section 1 (Corollary \ref{c} and Proposition \ref{p.4}) but the coalgebra is not right co-Frobenius. Also, this shows that left co-Frobenius does not imply neither uniqueness of right integrals nor existence of left integrals, since in this case, any combination of existence and uniqueness of integrals would imply the fact that $C$ is co-Frobenius by Theorem \ref{TheTh}.
\end{example}

\begin{example}\label{ex6}
Let $\Gamma$ be the directed graph (tree) obtained in the following way: start with the tree below $W$ (without a designated root):
$$
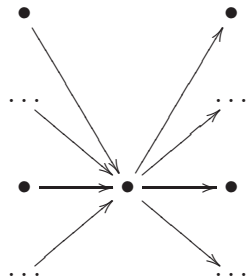
\xymatrix{
\bullet \ar[ddr] & & \bullet \\
\dots \ar[dr] & & \dots \\
\bullet\ar[r] & \bullet\ar[uur]\ar[ur]\ar[r]\ar[dr] & \bullet \\
\dots \ar[ur] & & \dots
}$$
This has infinitely many arrows going into the center-point $c$ and infinitely many going out. Then for each "free" vertex $x\neq c$ of this graph, glue (attach) another copy $W$ such that the vertex $x$ becomes the center of $W$, and one of the arrows of this copy of $W$ will be the original arrow $xc$ (or $cx$) with orientation. We continue this process for "free" vertices indefinitely to obtain the directed graph $\Gamma$ which has the property that each of its vertex $a$ has an infinite number of (direct) successors and an infinite number of predecessors.
Thus $|R(a)|=\infty$ and $L(a)=\infty$ for all the vertices $a$ of $\Gamma$, so we get $\dim(\int_{l,M})=\infty$ and $\dim(\int_{r,N})=\infty$ for all $M\in \Mm^C$ and $N\in {}^C\Mm$. Just as example \ref{ex1} this shows that the converse of Proposition \ref{p.3} does not hold even if we assume "existence" of left and right integrals; but the example here is non-cocommutative and has many types of isomorphism of simple comodules, and all spaces of integrals are infinite dimensional. 
\end{example}

Since integrals are tightly connected to the notions of co-Frobenius and QcF coalgebras, we also give some examples which show the fine non-symmetry of these notions; namely, we note that there are coalgebras which are QcF (both left and right), co-Frobenius on one side but not co-Frobenius. Also, it is possible for a coalgebra to be semiperfect (left and right) and QcF only on one side.

First, we note that the above general construction for graphs can be "enhanced" to contain non-pointed coalgebras. Namely, using the same notations as above, if $\Gamma$ is a labeled graph, i.e. a graph such that there is a positive natural number $n_v=n(v)$ attached to each vertex $v\in V$, then consider $K[\Gamma]$ to be the coalgebra with a basis $<(v_{ij})_{i,j=\overline{1,n(v)}}; (m_{ij})_{i=\overline{1,n_{l(m)}},j=\overline{1,n_{r(m)}}} \mid v\in V, m\in\Ee>$ and comultiplication and counit given by
\begin{eqnarray*}
\Delta(v_{ij}) & = & \sum_{k=1}^{n_v} v_{ik}\otimes v_{kj}\\
\Delta(m_{ij}) & = & \sum\limits_{k=1}^{n_{l(m)}}l(m)_{ik}\otimes m_{kj} + \sum_{k=1}^{n_{r(m)}}m_{ik}\otimes r(m)_{kj}\\
\varepsilon(v_{ij}) & = & \delta_{ij} \\
\varepsilon(m_{ij}) & = & 0
\end{eqnarray*}
Again, we can denote by $S_l(v,i)={}_K<v_{ki}\mid k=1,\dots,n_{v}>$ and $S_r(v,i)={}_K<v_{ik}\mid k=1,\dots,n_v>$; these will be simple left and respectively right $K[\Gamma]$-comodules. Also, let $E_l(v,i)={}_K<v_{ki}, k=1,\dots,n_v; m_{qi}, q=1,\dots,n_{l(m)}, m\in L(v)>$ and  put $E_r(v,i)={}_K<v_{ik}, k=1,\dots,n_v; m_{iq}, q=1,\dots,n_{r(m)}, m\in R(v)>$; these are the injective envelopes of $S_l(v,i)$ and $S_r(v,i)$ respectively. Let $S_{l/r}(v)=S_{l/r}(v,1)$ and $E_{l/r}(v)=E_{l/r}(v,1)$; these are representatives for the simple left/right $K[\Gamma]$-comodules, and for the indecomposable injective left/right $K[\Gamma]$-comodules.

\begin{example}\label{ex7}
Consider the labeled graph $\Gamma$ in the diagram bellow
$$\xymatrix{
\dots \ar[r] & {}_{(p_{-2})}\bullet^{a^{-2}} \ar[r]_{x^{-1}} & {}_{(p_{-1})}\bullet^{a^{-1}} \ar[r]_{x^{0}} & {}_{(p_0)}\bullet^{a^0} \ar[r]_{x^1} & {}_{(p_1)}\bullet^{a^1} \ar[r]_{x^2} & {}_{(p_2)}\bullet^{a^2} \ar[r]_{x^3} &  \dots
}$$
The vertices $a^n$ have labels positive natural numbers $p_n$ (they will be representing the simple subcoalgebras of the coalgebra $C=K[\Gamma]$ which are comatrix coalgebras of the respective size). Between each two vertices $a^{n-1},a^n$ there is a side $x^n$. The above coalgebra $C=K[\Gamma]$ then has a basis $\{a^n_{ij}\mid i,j=1,\dots,p_n,\,n\in\ZZ\}\sqcup \{x^n_{ij}\mid i=1,\dots,p_{n-1}, j=1,\dots,p_n,\, n\in\ZZ\}$ and structure
\begin{eqnarray*}
\Delta(a^n_{ij}) & = & \sum_{k=1}^{p_n}a^n_{ik}\otimes a^n_{kj} \\
\Delta(x^n_{ij}) & = & \sum_{k=1}^{p_{n-1}}a^{n-1}_{ik}\otimes x^n_{kj}+\sum_{k=1}^{p_n}x^n_{ik}\otimes a^n_{kj}\\
\varepsilon(a^n_{ij}) & = & \delta_{ij} \\
\varepsilon(x^n_j) & = & 0
\end{eqnarray*}
With the above notations, let $E_r(n)=E_r(a^{n})=E_r(a^{n},1)$, $E_l(n)=E_l(a^n)=E_l(a^n,1)$. We note that $E_l(n)^*\simeq E_r(n-1),\,\forall n$. First, note that if $M$ is a finite dimensional left $C$-comodule with comultiplication $\rho(m)=m_{-1}\otimes m_0$, $M^*$ is a right $C$-comodule with comultiplication $R$ such that $R(m^*)=m^*_0\otimes m^*_1$ if and only if 
\begin{eqnarray}
m_0^*(m)m^*_1 & = & m_{-1}m^*(m_0) \label{eqdual}
\end{eqnarray} 
This follows immediately by the definition of the left $C^*$-action on $M^*$. We then have the following formulas giving the comultiplication of $E_r(n-1)=<a^{n-1}_{1k}\mid 1\leq k\leq p_{n-1};\, x^n_{1k}\mid 1\leq k\leq p_n>$
\begin{eqnarray*}
a^{n-1}_{1k} & \mapsto & \sum_ja^{n-1}_{1j}\otimes a^{n-1}_{jk} \\
x^{n-1}_{1k} & \mapsto & \sum_ja^{n-1}_{1j}\otimes x^n_{jk}+\sum_jx^n_{1j}\otimes a^n_{jk} 
\end{eqnarray*}
and for $E_l(n)=<a^n_{k1}\mid 1\leq k\leq p_n; \, x^n_{k1}\mid 1\leq k\leq p_{n-1}>$ we have
\begin{eqnarray*}
a^n_{k1} & \mapsto & \sum_ja^n_{kj}\otimes a^n_{j1}\\
x^n_{k1} & \mapsto & \sum_ja^{n-1}_{kj}\otimes x^n_{j1}+\sum\limits_jx^n_{kj}\otimes a^n_{j1}
\end{eqnarray*}
Let $\{A^n_{k1}\mid 1\leq k\leq p_n;\,X^n_{k1}\mid 1\leq k\leq p_{n-1}\}$ be a dual basis for $E_l(n)^*$. Then, on this basis, the right comultiplication of $E_l(n)^*$ reads:
\begin{eqnarray*}
X^n_{k1} & \mapsto & \sum_i X^n_{i1}\otimes a^{n-1}_{ik}\\
A^n_{k1} & \mapsto & \sum_iX^n_{i1}\otimes x^n_{ik}+\sum_iA^n_{i1}\otimes a^n_{ik}
\end{eqnarray*}
Indeed, this can be easily observed by testing equation (\ref{eqdual}) for the dual bases $\{a^n_{k1}; x^n_{k1}\}$ and $\{A^n_{k1}; X^n_{k1}\}$. This shows that the 1-1 correspondence $a^{n-1}_{1k}\leftrightarrow X^n_{k1}$; $x^n_{1k}\leftrightarrow A^n_{k1}$ is an isomorphism of right $C$-comodules. Therefore, $E_l(n)^*\simeq E_r(n-1)$, and then also $E_r(n)\simeq E_l(n+1)^*$ for all $n$; thus we get that $C$ is QcF (left and right). One can also show this by first proving this coalgebra is Morita equivalent to the one obtained with the constant sequence $p_n=1$, which is QcF in a more obvious way, and using that "QcF" is a Morita invariant property. Note that $\dim(soc(E_r(n)))=\dim(<a^n_{1k}\mid 1\leq k\leq p_n>)=p_n$; $\dim(cosoc(E_r(n)))=\dim(soc(E_l(n+1)))=p_{n+1}$. Therefore, Remark \ref{1.rem} shows that $C$ is left co-Frobenius if and only if $(p_n)_n$ is an increasing sequence, and it is right co-Frobenius if and only if it is decreasing. Thus, taking $p_n$ to be increasing (decreasing) and non-constant we get a QcF coalgebra which is left (right) co-Frobenius and not right (left) co-Frobenius, and taking it to be neither increasing nor decreasing yields a QcF coalgebra which is not co-Frobenius on either side.
\end{example}

\begin{remark}
It is stated in \cite{Wm} (see also review MR1851217) that a coalgebra $C$ which is QcF on both sides must have left uniqueness of integrals ($\dim(\Hom^C(C,M))\leq \dim(M)$ for $M\in\Mm^C$). By Proposition \ref{p.3}, this is equivalent to $C$ being also left co-Frobenius. Nevertheless, by the above example we see that there are coalgebras which are both left and right QcF, but not co-Frobenius on either side. Note that even the hypothesis of $C$ being left QcF and right co-Frobenius would not be enough to imply the fact that left uniqueness of integrals holds. Some related uniqueness of integrals properties are stated in \cite{Wm2}; however, some questions arrise there as well, as noted also in Review MR2076973 (2005d:16070). \\
In fact, in the above example, denoting $S_{l/r}(n)=S_{l/r}(a_n,1)$, we have an exact sequence of right comodules $0\rightarrow S_r(n)\rightarrow E_r(n)\rightarrow S_r(n+1)\rightarrow 0$; also $K[\Gamma]^{K[\Gamma]}=\bigoplus\limits_nE_r(n)^{p_n}$ as right comodules. Therefore $\dim\Hom^{K[\Gamma]}(K[\Gamma],S_r(m))=\dim\prod\limits_n\Hom^{K[\Gamma]}(E_r(n),S_r(m))^{p_n}=p_{m-1}$ since $\Hom^{K[\Gamma]}(E_r(n),S_r(m))=\Hom^{K[\Gamma]}(S_r(n+1),S_r(m))=\delta_{n+1,m}$ (the $E_r(n)$'s are also local). Comparing this to $\dim(S_r(m))=p_m$, we see that any inequality is possible (unless some monotony property of $p_n$ is assumed, as above).
\end{remark}

\begin{example}\label{ex9}
Let $C=K[\Gamma]$ where $\Gamma=(V,\Ee)$ is the graph:
$$\xymatrix{
\bullet^{v_0} \ar[r]_{x_0} & \bullet^{v_1} \ar[r]_{x_1} & \bullet^{v_2} \ar[r]_{x_2} & \bullet^{v_3} \ar[r] & \dots
}$$
By the above considerations, $\dim(E_r(v_n))=2$ for all $n$ and $\dim(E_l(v_n))=2$ if $n>0$, $\dim(E_l(v_0))=1$. Also, $E_r(v_n)\simeq E_l(v_{n+1})^*$ for all $n$. These show that this coalgebra is semiperfect (left and right), it is left QcF (and even left co-Frobenius) but it is not right QcF since $E_l(v_0)$ is not isomorphic to any of the $E_r(v_n)^*$'s for any $n$ (by dimensions).
\end{example}

\begin{example}\label{ex10}
More generally, consider the labeled graph $\Gamma$
$$\xymatrix{
 {}_{(p_0)}\bullet^{a^0} \ar[r]_{x_0} & {}_{(p_1)}\bullet^{a^1} \ar[r]_{x_1} & {}_{(p_2)}\bullet^{a^2} \ar[r]_{x_2} & {}_{(p_3)}\bullet^{a^3} \ar[r]_{x_3} & \dots
}$$
and the corresponding coalgebra $C=K[\Gamma]$ (as constructed before example \ref{ex7}). In the same way as in examples \ref{ex7} and \ref{ex9} we get that $C$ is left QcF, but is not right QcF, and it is not left co-Frobenius if the sequence of natural numbers $(p_n)_n$ is not increasing (that is, if there is an $n$ such that $p_n>p_{n+1}$).
\end{example}

\begin{example}\label{ex11}
We note that there are coalgebras which are left and right semiperfect but are not left nor right QcF. Indeed, any finite dimensional coalgebra which is not QcF satisfies this (for example, the coalgebra of the graph $\bullet\rightarrow\bullet$).
\end{example}

The above examples also give a lot of information about the (left or right) semiperfect, QcF and co-Frobenius algebras: they show that any inclusion between such classes (such as, for example, left QcF and right semiperfect coalgebras into left and right QcF coalgebras) is a strict inclusion. Inclusions are in order since left (right) co-Frobenius coalgebras are left (right) QcF, and left (right) QcF coalgebras are left (right) semiperfect. We note however a very interesting fact:

\begin{theorem}\label{qcfsp}
Let $C$ be a left QcF coalgebra. Then $C$ is a semiperfect coalgebra (left and right).
\end{theorem}
\begin{proof}
We only need to check $C$ is right semiperfect. Let $S$ be a simple left comodule and assume $E(S)$ is infinite dimensional. Since it $C$ is left semiperfect, according to \cite[Proposition 2.3]{I}, the set $\{E(T)^*|T\in\Tt\}$ generates ${}^C\Mm$. Then, as each $E(T)^*$ is isomorphic to some $E(L)$, $L\in\Ss$ by Proposition \ref{p.1}, we get that there exists an epimorphism 
$$\bigoplus\limits_{i\in I}E(S_i)\stackrel{p}{\rightarrow}E(S)\rightarrow 0$$
where $S_i\in\Ss$ and $E(S_i)$ are finite dimensional and also $E(S_i)^*\cong E(T_i)$ for some $T_i\in\Tt$. Note that under our assumption there are infinitely many types of isomorphisms of simple comodules among $\{S_i\}_{i\in I}$, equivalently, there are infinitely many types of simples among $\{T_i\}_{i\in I}$. Indeed, assume there are only finitely many types of isomorphism of $S_i$'s. Since the $E(S_i)$'s are finite dimensional, there is a finite dimensional subcoalgebra $D$ of $C$ such that $\rho_{E(S_i)}(E(S_i))\subseteq E(S_i)\otimes D$ for all $i\in I$. Consequently, $D^\perp$ is contained in the annihilator of $\bigoplus\limits_{i\in I}E(S_i)$, so it is also in the annihilator of $P=\prod\limits_{i\in I}E(S_i)^*$. Since the annihilator of $P$ contains a closed ideal of finite codimension, 
$P$ is a rational $C^*$-module. It follows that $E(S)^*$ is rational too since it is a submodule of $P$. But $E(S)^*$ is also cyclic, and it would then be finite dimensional, a contradiction to the assumption.\\
Now fix some $i\in I$. We may obviously assume that the kernel of $p$ contains none of the $E(S_i)$'s (see also \cite[Proposition 2.4]{I} for a more general statement). This means that $p(E(S_i))\neq 0$, so $S\subseteq p(E(S_i))$, thus $S$ appears as a factor in the Jordan-Holder series of $E(S_i)$. \\
We can choose a subcomodule $X_i\subseteq E(S_i)$ of minimal dimension such that $X_i$ projects onto $S$, equivalently, $S$ occurs on top on a Jordan-Holder series of $X_i$. Thus there is a subcomodule $K_i\subset X_i$ such that $X_i/K_i\simeq S$. We see that $X_i$ is local. Indeed, if $M$ is another maximal subcomodule of $X_i$ with $M\neq K_i$, then as $X_i/K_i=M+K_i/K_i\cong M/M\cap K_i$, we get that $M$ projects onto $S$ but has dimension less than ${\rm dim}(X_i)$. Therefore, dualy, we obtain that the socle of $X_i^*$ is simple isomorphic to $S^*$, so there is a monomorphism $X_i^*\hookrightarrow E(S^*)$. Also, since the socle of $X_i\subseteq E(S_i)$ is equal to $S_i$, we see that $S_i^*$ occurs in the Jordan-Holder series of $X_i^*$, and so also in a Jordan-Holder series of $E(S^*)$. But since there are infinitely many types of isomorphisms of simples among the $S_i$'s, we get that $E(S^*)$ is infinite dimensional, which contradicts the fact that $C$ is left semiperfect. This ends our proof. 
\end{proof}



The following table sums up all these examples with respect to the left or right semiperfect, QcF or co-Frobenius coalgebras. More precisely, one has \\
$\{$Left (Right) co-Frobenius$\}\,\subset\,\{$Left (Right) QcF$\}\,\subseteq\,\{$Left and Right semiperfect$\}$, \\
and all inclusions between suitable left and right combinations of these are strict. In the column to the left, the coalgebras in the various examples are refered (Example \ref{ex7} contains examples for more than one situation); the other columns record the properties of these coalgebras. Moreover, all these coalgebras are non-cocommutative. Further minor details here are left to the reader.


\begin{center}
\begin{tabular}{||c|c|c|c|c|c|c||}\hline 
Example & left & right & left QcF & right QcF & left & right \\ 
& semiperfect & semiperfect & & & co-Frob & co-Frob \\ \hline
Ex \ref{ex2.5} & $\surd$ & $\surd$ & & $\surd$ & & $\surd$ \\ \hline
Ex \ref{ex3}   & & $\surd$ & & & & \\ \hline
Ex \ref{ex6}   & & & & & & \\ \hline
Ex \ref{ex7}(1) & $\surd$ & $\surd$ & $\surd$ & $\surd$ & $\surd$ & \\ \hline
Ex \ref{ex7}(2) & $\surd$ & $\surd$ & $\surd$ & $\surd$ & & \\ \hline
Ex \ref{ex7}(3) & $\surd$ & $\surd$ & $\surd$ & $\surd$ & $\surd$ & $\surd$\\ \hline
Ex \ref{ex9}   & $\surd$ & $\surd$ & $\surd$ & & $\surd$ & \\ \hline
Ex\ref{ex10} & $\surd$ & $\surd$ & $\surd$ & & & \\ \hline
Ex\ref{ex11} & $\surd$ & $\surd$ & & & & \\ \hline
\end{tabular}
\end{center}


\section{Further Applications to Categorical Characterizations \\ and Integrals for Algebras}

A QF algebra is characterized by the property that injective left (or right) $A$-modules coincide with the projective ones. In fact, for $A$ to be Quasi-Frobenius it is enough to ask that any injective is projective, or that any projective is injective (Faith). For coalgebras, a left QcF coalgebra is characterized by the fact that any injective is projective (see e.g. \cite{DNR}), but this is not equivalent to the fact that projectives are injective, since it may happen that the only projectives are $0$. This is the case for example for the divided power coalgebra $C=K\{c_n|n\in\NN\}$, $\Delta(c_n)=\sum\limits_{i+j=n}c_i\otimes c_j$, $\varepsilon(c_n)=\delta_{0,n}$. Any comodule over this coalgebra is a direct sum of indecomposable comodules isomorphic either to $C$ or to the $n$'th term $C_n$ of the coradical filtration of $C$. None of these is projective. Also, the condition ``injective$\rightarrow$projective'' is not enough for two-sided QcF. Hence for the characterization of the two-sided Quasi-co-Frobenius coalgebras it is to be expected that a symmetric such condition should be required. Theorem \ref{qcfsp} indeed allows us to prove:

\begin{theorem} \label{ThPrIn}
The following conditions are equivalent for a coalgebra $C$:\\
(i) $C$ is QcF.\\
(ii) The class of projective left $C$-comodules coincides with the class of injective left $C$-comodules.\\
(iii) Projective right $C$-comodules coincide with injective right $C$-comodules.   
\end{theorem}
\begin{proof}
(i)$\Rightarrow$(ii) We only need to show that projectives are injective. Let $P\in{}^C\Mm$ be projective. Let $C^{(I)}\rightarrow P\rightarrow 0$ be an epimorpism: it exists, since in this case ${}^CC$ is a generator (e.g. by \cite{NT1}). By projectivity, $P$ is a direct summand in $C^{(I)}$ which is an injective comodule, so $P$ is injective.\\
(ii)$\Rightarrow$(i) Since ${}^CC$ is projective, it follows that $C$ is right QcF. By Theorem \ref{qcfsp}, $C$ is then left semiperfect too. Thus, the injective envelopes $E(T)$ of right simple comodules $T$ are finite dimensional. Hence, $E(T)^*$ is a projective left $C$-comodule for any $T\in\Tt$, so it is injective. Thus, since it is indecomposable, it follows that $E(T)^*\cong E(S)$ for some $S\in\Ss$, and so, by Proposition \ref{p.1} it follows that $C$ is left QcF too. 
\end{proof}

\subsection*{Semiperfect coalgebras}

It seems worthwhile at this point to note a functorial-categorical characterization of semiperfect coalgebras which paralels those of (quasi)-co-Frobenius coalgebras by integrals. We note that for a co-Frobenius coalgebra, the forgetful functor $f.d.{}^C\Mm\rightarrow {}_K\Mm$ is the restriction of a representable functor $\int_{r,-}:{}^C\Mm\rightarrow {}_K\Mm$. We first need a small lemma:

\begin{lemma}
Let $C$ be a left semiperfect coalgebra. Then $$Rat(\prod\limits_{T\in\Tt}E(T)^{p(T)})=\bigoplus\limits_{T\in\Tt}E(T)^{p(T)}=C$$
\end{lemma}
\begin{proof}
Let $x=(x_T)_{T\in\Tt}\in Rat(\prod\limits_{T\in\Tt}E(T)^{p(T)})$; then the annihilator of $x$ contains a closed ideal of finite codimension $D^\perp$, $D$ finite dimensional subcoalgebra of $C$. Then $D^\perp x=0$, so $D^\perp\cdot x_T=0$ for all $T$ and so $D^\perp \cdot C^*x_T=0$. Since $C^*x_T$ are rational finite dimensional left $C^*$-submodules of $E(T)^{p(T)}$, whenever $x_T\neq 0$, we have that the socle of $C^*\cdot x_T$ contains a simple right $C$-comodule isomorphic to $T$. Therefore, in this case, it follows that $D^\perp$ cancels $T$ i.e. $D^\perp\cdot T=0$, and so, $T\subseteq D$. This shows that only finitely many $x_T$'s are nonzero, and this ends the proof since the converse inclusion is obvious.
\end{proof}

\begin{theorem}
Let $C$ be a coalgebra. Then $C$ is left semiperfect if and only if the forgetful functor $f.d.{}^C\Mm\rightarrow {}_K\Mm$ is the restriction of a representable functor ${}^C\Mm\rightarrow {}_K\Mm$.
\end{theorem}
\begin{proof}
If $C$ is semiperfect, let $P=\bigoplus\limits_{T\in\Tt}E(T)^*{}^{p(T)}$. Then we have
$$Rat(P^*)=Rat((\bigoplus\limits_{T\in\Tt}(E(T)^*)^{p(T)})^*)=Rat(\prod\limits_{T\in\Tt}(E(T)^*)^{*p(T)})=Rat(\prod\limits_{T\in\Tt}E(T)^{p(T)})=C$$ because of the previous Lemma. Thus we have the natural isomorphisms for finite dimensional left comodules $N$
\begin{eqnarray*}
\Hom^C(P,N) & \cong & \Hom^C(P,(N^*)^*)\cong\Hom_{C^*}(N^*,P^*)\cong \\
& \cong & \Hom_{C^*}(N^*,Rat(P^*))\cong\Hom^C(N^*,C)\cong(N^*)^*=N.
\end{eqnarray*}
Conversely, assume $\Hom^C(P,N)\cong N$ naturally as vector spaces for all finite dimensional left comodules $N$. Then for right comodules $M$ we have natural isomorphisms:
\begin{eqnarray}
\Hom^C(M,Rat(P^*))\cong\Hom_{C^*}(M,P^*)=\Hom_{C^*}(P,M^*)\cong M^*\cong \Hom^C(M,C) \label{eq3}
\end{eqnarray}
Since $\Hom(-,Rat(P^*))$ is exact on sequences of finite dimensional comodules, it follows that $Rat(P^*)$ is an injective right $C$-comodule (see \cite[Section 2.4]{DNR}). Then for a right comodule $M$, write $M=\lim\limits_{\stackrel{\rightarrow}{i}}M_i$ as a directed limit of its finite dimensional subcomodules $M_i$. Then there is a natural isomorphism
\begin{eqnarray*}
\Hom(M,Rat(P^*)) & = & \Hom(\lim\limits_{\stackrel{\rightarrow}{i}}M_i,Rat(P^*)) \\
& = & \lim\limits_{\stackrel{\leftarrow}{i}}\Hom(M_i,Rat(P^*)) \,\,\,(Rat(P^*){\rm\,injective}) \\
& = & \lim\limits_{\stackrel{\leftarrow}{i}}\Hom(M_i,C) \,\,\,({\rm\,by\,the\,naturality\,of\,the\,isomorphism\,in\,}(\ref{eq3}))\\
& = & \Hom(\lim\limits_{\stackrel{\rightarrow}{i}}M_i,C)=\Hom(M,C) \,\,\,(C{\rm\,injective})
\end{eqnarray*}
It then follows that $C\cong Rat(P^*)$ as right comodules. In particular, the monomorphism $C\hookrightarrow P^*$ of left $C^*$-modules gives rise to a morphism of right $C^*$-modules $\psi:P\rightarrow C^*$ (as in Lemma \ref{l.1}) which has dense image. Since $\psi(P)\subseteq Rat(C^*_{C^*})$, this shows that $Rat(C^*_{C^*})$ is dense in $C^*$ and so $C$ is left semiperfect (see e.g. \cite{L} or \cite[Section 3.2]{DNR}).
\end{proof}

\begin{remark}
Note that this is the best result possible: we cannot ask that the forgetful functor is actually representable. If $\Hom^C(H,-)=U:f.d.{}^C\Mm\rightarrow {}_K\Mm$ with $H\in f.d.{}^C\Mm$, then when $C$ is left semiperfect infinite dimensional, then there are infinitely many types of simple left comodules. But then since $H$ is finite dimensional, there is some simple $S$ left comodule such that $\Hom^C(H,S)=0$ (for example, an $S$ not in the Jordan-Holder series of $H$), so $\Hom^C(H,S)\not\cong S$.
\end{remark}

\subsection*{Integrals for algebras}

If $A$ is a topological algebra with a linear topology which has a basis of neighbourhoods of $0$ consisting of ideals of finite codimension, we call such an algebra a topological algebra of algebraic type, or AT-algebra (see also \cite{II}). This is motivated by the situation when one is interested in in the study of only a certain subcategory $\Cc$ of that of finite dimensional representations of $A$. If this the case, we can introduce the linear topology $\gamma$ on $A$ with a basis of open neighbourhoods of $0$ consisting of all ideals $I$ which are the annihilator of some $M\in\Cc$. If $\Cc$ is closed under subobjects, quotients and direct sums, then $\Cc$ coincides with the category of finite dimensional topological $A$-modules. Moreover, for an AT-algebra, the category of finite dimensional topological left $A$-modules is equivalent to the cateogory of finite dimensional right $R(A)$-comodules, where $R(A)$ is the coalgebra of representative functions on $A$, that is, $R(A)=\{f:A\rightarrow K | \Ker(f) {\rm\,contains\,a\,cofinite\,closed\,ideal\,}\}$, equivalently, $R(A)=\{f:A\rightarrow K{\rm\,continuous}\,|\, \exists g_i,h_i{\rm\,continuous},\,i=1,\dots,n,\,{\rm with\,}f(xy)=\sum\limits_{i=1}^ng_i(x)h_i(y),\,\forall x,y\in A\}$, where $K$ is considered with the discrete topology. Furthermore, $R(A)$ is the span of all functions $\eta_{ij}$ where $\eta:A\rightarrow \End(V)$ is a continuous finite dimensional representation of $A$, $v_i$ is a basis of $V$ and $\eta(a)\cdot v_i=\sum\limits_{j}\eta_{ji}(a)v_j$. Hence $R(A)=\lim\limits_{\stackrel{\longleftarrow}{I{\rm\,open\,ideal}}}(\frac{A}{I})^*$ as a coalgebra (see also \cite[Section 2.5]{DNR}). For such an algebra, we can define an ``integral'' functor, i.e. an integral space for each finite dimensional topological left $A$-module $M$ by $\int_{l,M}=\Hom_A(R(A),M)$. For an arbitrary algebra $A$ (with no specified such topology) we can consider the topology having all cofinite ideals as a basis of neighbourhoods of $0$, and then $R(A)=A^0$ the finite dual coalgebra. The definition then still makes sense for this case, and all the results on integrals apply to these situations. One would then be entitled to call an AT-algebra $A$ weakly (quasi)-Frobenius if $R(A)$ is a (quasi)-co-Frobenius coalgebra (see \cite{II}). Hence, for example, among other, we would have:

\begin{corollary}
Let $A$ be a topological algebra of algebraic type, and let $F({}_A\Mm)$ ($F(\Mm_{A})$) denote the category of left (right) $A$-modules which are the sum of their finite dimensional submodules, i.e. modules which are directed limits of finite dimensional $A$-modules. (then $F({}_A\Mm)=Rat({}_{R(A)^*}\Mm)$). Then the following are equivalent. \\
(i) $A$ is weakly quasi-Frobenius.\\
(ii) In $F({}_A\Mm)$ injectives and projectives coincide. \\
(iii) In $F(\Mm_A)$ injectives and projectives coincide.
\end{corollary}

\begin{corollary}
The following are equivalent for a topological algebra $A$ of algebraic type.\\
(i) A is weakly Frobenius.\\
(ii) $\dim(\int_{l,M})=\dim(M)$ for all left finite dimensional left $A$-modules $M$.\\
(iii) $\dim(\int_{r,M})=\dim(M)$ for all left finite dimensional right $A$-modules $M$.
\end{corollary}

\begin{corollary}[of Theorem \ref{qcfsp}]
Let $A$ be an (AT-) algebra. If the category of (topological) finite dimensional left $A$-modules has enough injectives (respectively projectives) and every injective is projective (respectively, every projective is injective), then the category of (topological) finite dimensional right $A$-modules has enough injectives (respectively, projectives) too. 
\end{corollary}

\section{Locally Compact Groups}

We examine what do the generalized integrals represent for the case of locally compact groups. First we look at a very simple case. Consider the measure $d\mu_t(x)=e^{itx}dx$ on the group $(\RR,+)$ for some $t\in \RR$, that is, $\int_\RR f(x)d\mu(x)=\int_\RR f(x)e^{itx}dx$ for any $f\in L^1(\RR)$. Then this measure $\mu_t$ has a special type of "invariance", since $\int_\RR f(x+a)d\mu_t(x)=\int_\RR f(x+a)e^{itx}dx=\int_\RR(f(x)e^{it(x-a)})dx=e^{-ita}\int_\RR f(x)d\mu_t(x)$. Equivalently, this means that for any Borel set $U$ we have $\mu_t(U+a)=\int_\RR \chi_U(x-a)d\mu_t(x)=e^{ita}\mu_t(U)$, that is, translation by $a$ of a set has the effect of "scaling" it (its measure) by $e^{ita}$ (note that here $t$ could be any complex number). We generalize this for a locally compact group $G$ with left invariant Haar measure $\lambda$. 

Let $\mu$ be a complex vector measure on $G$, that is, $\mu=(\mu_1,\dots,\mu_n)$ and so $\mu(U)=(\mu_1(U),\dots,\mu_n(U))\in \CC^n$ for any Borel subset $U$ of $G$. We will be looking at the above type of invariance for such a measure $\mu$, that is, measures such that right translation of $U$ by $g\in G$ will have the effect of scaling $\mu(U)$ by $\eta(g)$, where $\eta(g)$ must be an $n\times n$ matrix, i.e. $\mu(U\cdot g)=\alpha(g)\cdot \mu(U)$. With the natural left action of $G$ on the set of all functions $f:G\rightarrow \CC$ defined by $(y\cdot f)(x)=f(xy)$, this writes equivalently
\begin{eqnarray*}
\int\limits_Gg\cdot\chi_Ud\mu& = & \int\limits_G\chi_U(x\cdot g)d\mu(x)=\int\limits_G\chi_{Ug^{-1}}d\mu\\
& = & \mu(U\cdot g^{-1})=\alpha(g^{-1})\mu(U) \\
& = & \alpha(g^{-1})\int\limits_G\chi_Ud\mu
\end{eqnarray*}
This is extended to all $L^1$ functions $f$ by $\int_Gx\cdot fd\mu=\eta(x)\int_Gfd\mu$, where $\eta(x)=\alpha(x^{-1})$. Note that we have $\eta(xy)\int_Gfd\mu=\int_Gxy\cdot fd\mu=\eta(x)\int_Gy\cdot fd\mu=\eta(x)\eta(y)\int_Gfd\mu$. This leads to the following definition:

\begin{definition}
Let $G$ be a topological group and $\int:C_c(G)\rightarrow V=\CC^n$ be a linear map, where $C_c(G)$ is the space of continuous functions of compact support on $G$. We say $\int$ is a quantum $\eta$-invariant integral (quantified by $\eta$) if $\int (x\cdot f)=\eta(x)\int(f)$ for all $x\in G$, where $\eta:G\rightarrow \End(V)$.
\end{definition}

We note that the "quantum" factor $\eta$ is a representation which is continuous if the linear map $\int$ is itself continuous, where the topology on $C_c(G)$ is that of uniform convergence, thus that given by the sup norm $||f||=\sup\limits_{x\in G}|f(x)|$ for $f\in C_c(G)$. For example, by general facts of measure theory, if $G$ is locally compact $\int=(\lambda_1,\dots,\lambda_n)$ is continuous whenever $\int=\int d\mu$, $\mu=(\mu_1,\dots,\mu_n)$ and $\mu_i$ are positive measures (i.e. $\lambda_i=\int(-)d\mu_i$ is positive in the sense that $\lambda_i(f)\geq 0$ whenever $f\geq 0$) or more generally, $\mu_i=\mu_{i1}-\mu_{i2}+i(\mu_{i3}-\mu_{i4})$, where $\mu_{ij}$ are all positive measures, since any $\sigma$-additive complex measure $\mu_i$ can be written like this. Similarly, it can be seen that any continuous $\lambda_i:C_c(G)\rightarrow \CC$ can be written as $\lambda_i=(\lambda_{i1}-\lambda_{i2})+i(\lambda_{i3}-\lambda_{i4})$ with $\lambda_{ij}$ positive linear functionals which can be represented as $\lambda_{ij}=\int(-)d\mu_{ij}$ by the Riesz Representation Theorem. We refer to \cite{Ru} for the facts of basic measure theory.

\begin{proposition}
Let $G$ be a locally compact group, $\int:C_c(G)\rightarrow V=\CC^n$ be a quantum $\eta$-invariant integral and $W=\im(\int)$. Then:\\
(i) $W$ is an $\eta$-invariant subspace of $V$, that is, $\eta(x)W\subseteq W$ for all $x\in G$.\\
Consider the induced map $\eta:G\rightarrow \End(W)$ which can be considered by (i), and denote it also by $\eta$. Then:\\
(ii) $\eta$ is a representation of $G$, so $\eta:G\rightarrow GL(W)$.\\
(iii) $\eta$ is a continuous representation if $\int$ is continuous.
\end{proposition}
\begin{proof}
(i) For $w=\int(f)\in W$, $\eta(x)w=\eta(x)\int(f)=\int(x\cdot f)\in W$.\\
(ii) For any $w=\int(f)\in W$ and $x,y\in G$ we have $\eta(xy)w=\eta(xy)\int(f)=\int(xy\cdot f)=\eta(x)\int(y\cdot f)=\eta(x)\eta(y)\int(f)=\eta(x)\eta(y)w$ and so $\eta(xy)=\eta(x)\eta(y)$ (since here $\eta$ is considered with values in $\End(W)$). Since $1\cdot f=f$, we get $\eta(1)={\rm Id}_W$. Hence ${\rm Id}_W=\eta(1)=\eta(xx^{-1})=\eta(x)\eta(x^{-1})=\eta(x^{-1}\eta(x))$ so $\eta:G\rightarrow GL(W)$ is a representation. We note that this result holds also for the case when $V$ is an infinite dimensional complex vector space.\\
(iii) Let $w_1,\dots,w_k$ be an orthonormal basis of $W$ and let $\varepsilon$ be fixed. For each $i$, let $f_i\in C_c(G)$ be such that $w_i=\int f_i$. Since $\int$ is continuous, we can choose $\delta_i$ such that $|\int(g-f_i)|<\varepsilon$ whenever $||g-f_i||<\delta_i$, and let $\delta=\min\{\delta_i\mid i=1,\dots,n\}/2$. Since $f_i$ is compactly supported and continuous, it is (well known that it is) also uniformly continuous, and therefore there is a neighbourhood $U_i$ of $1$ - which may be assumed symmetric - such that if $y^{-1}z\in V_i$ then $|f_i(z)-f_i(y)|<\delta$. Therefore, if $x\in U_i$ and $y\in G$, $|(x\cdot f_i-f_i)(y)|=|f_i(yx)-f_i(y)|<\delta$ so $||x\cdot f_i-f_i||\leq \delta<\delta_i$. Hence, 
$$|(\eta(x)-{\rm Id})(w_i)|=|\eta(x)\int(f_i)-\int(f_i)|=|\int(xf_i-f_i)|<\varepsilon$$
Then this holds for all $x\in U=\bigcap\limits_{i=1}^kU_i$. For any $w\in W$, $w=\sum\limits_{i=1}^ka_iw_i$ and for all $x\in U$ we have
\begin{eqnarray*}
|(\eta(x)-{\rm Id})(w)| & = & |\sum\limits_{i=1}^ka_i(\eta(x)-{\rm Id})(w_i)| \leq \sum\limits_{i=1}^k|a_i|\cdot|(\eta(x)-{\rm Id})(w_i)|\\
& \leq & \varepsilon\sum\limits_{i=1}^k|a_i|\leq \varepsilon\cdot \sqrt{k\sum\limits_{i=1}^ka_i^2} \leq \varepsilon\sqrt{k}||w||
\end{eqnarray*}
This shows that the norm of $\eta(x)-{\rm Id}$ (as a continuous linear operator on $W$) is small enough for $x\in U$: $||\eta(x)-\eta(1)||=||\eta(x)-{\rm Id}||\leq \varepsilon\sqrt{k}$, so $\eta$ is continuous at $1$, so it is continuous everywhere since it is a morphism of groups $\eta:G\rightarrow GL(W)$.
\end{proof}

Now let $C=R_c(G)$ be the coalgebra (and actually Hopf algebra) of continuous representative functions on $G$. It is well known that any continuous (not necessary unitary) finite dimensional representation of $\eta:G\rightarrow GL(V)\subset\End(V)$ becomes a right $R_c(G)$-comodule in the following way: if $(v_i)_{i=1,n}$ is a basis of $V$, one writes 
\begin{eqnarray}\label{e1}
g\cdot v_i=\sum\limits_{j}\eta_{ji}(g)v_j 
\end{eqnarray}
and then it is straightforward to see that 
$\eta_{ij}(gh)=\sum\limits_{k}\eta_{ik}(g)\eta_{kj}(h)$
so $\eta_{ij}\in R_c(G)$ and $\rho:V\rightarrow V\otimes R_c(G)$, 
\begin{equation}\label{e2}
v_i\mapsto \sum\limits_{j}v_j\otimes \eta_{ij}
\end{equation} 
is a comultiplication. Conversely, the action of (\ref{e1}) defines a representation of $G$ whenever $V$ is a finite dimensional $R_c(G)$-comodule defined by (\ref{e2}). Also, the formula in (\ref{e1}) defines a continuous representation, since the linear operations on $V$ - a complex vector space with the usual topology - are continuous, and the maps $\eta_{ij}$ are continuous too. Moreover, $\varphi:V\rightarrow W$ is a (continuous) morphism of left $G$-modules (representations) if and only if it is a (continuous) morphism of right $R_c(G)$-comodules, that is, the categories of finite dimensional right $R_c(G)$-comodules and of finite dimensional $G$-representations are equivalent. We can now give the interpretation of generalized algebraic integrals for locally compact groups:

\begin{proposition}
Let $\eta:G\rightarrow \End(V)$ be a (continuous) finite dimensional representation of $G$ and $\int:C_c(G)\rightarrow V=\CC^n$ be an $\eta$-invariant integral as above, i.e. $$\int x\cdot f=\eta(x)\int f$$ and let $\lambda:R_c(G)\rightarrow V$ be the restriction of $\int$ to $R_c(G)\subseteq C_c(G)$: $\lambda(f)=\int f$. Then $\lambda\in\int_{V}=\Hom^{R_c(G)}(R_c(G),V)$ (in the sense of Definition \ref{def}).
\end{proposition}
\begin{proof}
It is enough to show that $\lambda$ is a morphism of left $G$-modules. But this is true, since $x\cdot\lambda(f)=\eta(x)\int\limits_Gfd\mu=\int\limits_Gx\cdot fd\mu=\lambda(x\cdot f)$. 
\end{proof}

We finish with a theorem for uniqueness and existence of $\eta$-invariant integrals. First we note that, as an application of a purely algebraic result, we can get the following nice and well known fact in the theory of compact groups:

\begin{proposition}
Let $G$ be a compact group. Then any finite dimensional continuous (not necessary unitary) representation $\eta:G\rightarrow GL(V)$ of $G$ is completely reducible.
\end{proposition}
\begin{proof}
By the above comments, the statement is equivalent to showing that $V$ is cosemisimple as a $R_c(G)$-comodule. But $R_c(G)$ is a Hopf algebra $H$ whose antipode $S$ has $S^2={\rm Id}$ (since $S(f)(x)=f(x^{-1})$) and it has integrals (in the Hopf algebra sense) - the left Haar integral, as it also follows by the above Proposition. This integral is nonzero and defined on all $f\in R_c(G)$, since $G$ is compact. Then a result of \cite{Su} (with a very short proof also given in \cite{DNT}) applies, which says that an involutory Hopf algebra with non-zero integrals is cosemisimple. Therefore, $R_c(G)$ is cosemisimple so $V$ is cosemisimple.
\end{proof}

\begin{remark}
It is well known that any continuous representation $V$ of $G$ is completely reducible, but for infinite dimensional representations, the decomposition is in the sense of Hilbert direct sums of Hilbert spaces.
\end{remark}

For a continuous finite dimensional representation $\eta:G\rightarrow GL(V)$, let ${\rm Int}_\eta$ denote the space of all continuous quantum $\eta$-invariant integrals on $C(G)$. Then we have:

\begin{theorem}[Uniqueness of quantum invariant integrals] \,$\underline{}$\\ 
Let $G$ be a compact group, and $\eta:G\rightarrow GL(V)$ a (continuous) representation of $G$. Then 
$$\dim({\rm Int}_\eta)\leq\dim V$$
\end{theorem}
\begin{proof}
By the Peter-Weyl theorem, it is known that the continuous representative functions $R_c(G)$ are dense in the space of all continuous functions $C(G)$ in the uniform norm. Therefore, the morphism of vector spaces ${\rm Int}_\eta\rightarrow \int_{l,V}$ given by the restriction is injective. Since $\dim(\int_{l,V})=\dim(V)$ by Theorem \ref{TheTh}, we get the conclusion. 
\end{proof}

\begin{remark}
In particular, we can conclude the uniqueness of the Haar integral in this way. However, the existence cannot be deduced from here, since, while the uniqueness of the Haar measure is not an essential feature of this part of the Peter-Weyl theorem, the existence of the left invariant Haar measure on $G$ is. 
\end{remark}

The existence of quantum integrals can be easily obtained constructively from the existence of the Haar measure as follows. We note that for any $v\in V$, we can define $H_\eta(v)\in{\rm Int}_\eta$ by $H_\eta(v)(f)=\int\limits_Gf(x^{-1})\eta(x)\cdot v$, where $\int\limits_G$ is the (a) left translation invariant Haar measure on $G$. Let $e_1,\dots,e_n$ be a $\CC$-basis of $V$ and $\eta(x)=(\eta_{i,j}(x))_{i,j}$ be the coordinates formula for $\eta$. Then $H_\eta(v)(f)=(\sum\limits_{j=1}^n\int\limits_Gf(x^{-1})\eta_{i,j}(x)v_j)_{i=1}^n$. We note that this is well defined, that is, $H_\eta(v)$ is a quantum $\eta$-invariant integral: indeed let $\Lambda=H_\eta(v)=(\Lambda_i)_{i=1}^n$; then using the fact that $\int\limits_G$ is the left translation invariant Haar integral, we get
\begin{eqnarray*}
\lambda_i(x\cdot f) & = & \sum\limits_{j=1}^n\int\limits_Gf(y^{-1}x)\eta_{ij}(y)v_jdy \,\,\, ({\rm substitute}\,y=xz)\\
& = & \sum\limits_{j=1}^n\int\limits_G f(z^{-1})\eta_{ij}(xz)v_jdz \,\,\,(\int\limits_G{\rm \,is\,the\,left\,Haar\,integral})\\
& = & \sum\limits_{k=1}^n\sum\limits_{j=1}^nf(x^{-1})\eta_{ik}(x)\eta_{kj}(z)v_jdz \,\,\,(\eta{\rm\,is\,a\,morphism})\\
& = & \sum\limits_{k=1}^n\eta_{ik}(x)\left(\sum\limits_{j=1}^n\int\limits_G f(z^{-1})\eta_{kj}(z)v_jdz \right)\\
& = & \sum\limits_{k=1}^n\eta_{ik}(x)\Lambda_k(f)
\end{eqnarray*}
so $\Lambda(x\cdot f)=\eta(x)\cdot \Lambda(f)$, i.e. $\Lambda\in{\rm Int}_\eta$. Define $\theta_\eta:{\rm Int}_\eta\rightarrow V$ by $\theta_\eta(\Lambda)=(\sum\limits_{i=1}^n)\Lambda_j(\eta_{ij})_{i=1}^n$. Then we have $\theta_\eta\circ H_\eta={\rm Id_V}$. Indeed, we have
\begin{eqnarray*}
\theta_\eta(H_\eta(v)) & = & \sum\limits_{j=1}^nH_\eta(v)_j(\eta_{ij})\\
& = & \sum\limits_{j=1}^n\sum\limits_{k=1}^n\int\limits_G\eta_{ij}(x^{-1})\eta_{jk}(x)v_kdx\\
& = & \sum\limits_{k=1}^n\int\limits_G\eta_{ik}(1_G)v_kdx \,\,\,(\eta(x^{-1})\eta(x)=\eta(1_G)={\rm Id}_n\in\mathcal{M}_n(\CC))\\
& = & \sum\limits_{k=1}^n\int\limits_G\delta_{ik}v_kdx=\int\limits_Gv_idx=v_i
\end{eqnarray*} 
where we have assumed, without loss of generality, that $\int\limits_G$ is a normalized Haar integral.

\bigskip\bigskip\bigskip\bigskip



\vspace*{3mm} 
\begin{flushright}
\begin{minipage}{148mm}\sc\footnotesize

Miodrag Cristian Iovanov\\
University of Southern California\\
Department of Mathematics, 3620 South Vermont Ave. KAP 108 \\
Los Angeles, California 90089-2532 \&\\
University of Bucharest, Faculty of Mathematics, Str.
Academiei 14,
RO-70109, Bucharest, Romania
{\it E--mail address}: {\tt
yovanov@gmail.com; iovanov@usc.edu}\vspace*{3mm}

\end{minipage}
\end{flushright}
\end{document}